\documentclass[11pt]{amsart}%
\usepackage{graphicx}
\usepackage{amscd}
\usepackage{amsmath}
\usepackage{mathrsfs}
\usepackage{mathtools}
\usepackage{amsfonts}
\usepackage{enumerate}
\usepackage{amssymb}
\usepackage{color}%
\providecommand{\U}[1]{\protect\rule{.1in}{.1in}}
\newtheorem{theorem}{Theorem}
\theoremstyle{plain}

\newtheorem{lemma}{Lemma}[section]

\newtheorem{remark}{Remark}[section]

\numberwithin{equation}{section}
\numberwithin{theorem}{section}

\everymath{\displaystyle}
\begin{document}
\date{July 14, 2014 (File: XmLaguerreTypesFV.tex)}
\title[Three Types of Exceptional $X_{m}$-Laguerre Polynomials]{A New Class of Exceptional Orthogonal Polynomials: The Type III $X_{m}%
$-Laguerre Polynomials And The Spectral Analysis of Three Types of Exceptional
Laguerre Polynomials}
\author{Constanze Liaw}
\address{Department of Mathematics, Baylor University, One Bear Place \#97328, Waco, TX 76798-7328}
\email{Constanze\_Liaw@baylor.edu}
\author{Lance L.~Littlejohn}
\address{Department of Mathematics, Baylor University, One Bear Place \#97328, Waco, TX 76798-7328}
\email{Lance\_Littlejohn@baylor.edu}
\author{Robert Milson}
\address{Department of Mathematics and Statistics, Dalhousie University\\
Halifax, Nova Scotia, Canada}
\email{milson@mathstat.dal.ca}
\author{Jessica Stewart}
\address{Department of Mathematics, Baylor University, One Bear Place \#97328, Waco, TX 76798-7328}
\email{Jessica\_Stewart@baylor.edu}
\thanks{Constanze Liaw is partially supported by the NSF grant DMS-1261687.}
\subjclass[2010]{Primary: 33C45, 34B24, 42C05. Secondary: 33C47, 47B25.}
\keywords{exceptional orthogonal polynomials, Laguerre polynomials, complete
eigenfunctions, second-order ordinary differential expression, spectral
theory, self-adjoint operator, root asymptotics.}

\begin{abstract}
The Bochner Classification Theorem (1929) characterizes the polynomial
sequences $\{p_{n}\}_{n=0}^{\infty}$, with $\deg p_{n}=n$ that simultaneously
form a complete set of eigenstates for a second-order differential operator
and are orthogonal with respect to a positive Borel measure having finite
moments of all orders. Indeed, up to a complex linear change of variable, only
the classical Hermite, Laguerre, and Jacobi polynomials, with certain
restrictions on the polynomial parameters, satisfy these conditions. In 2009,
G\'{o}mez-Ullate, Kamran, and Milson found that for sequences $\{p_{n}%
\}_{n=1}^{\infty}$, $\deg p_{n}=n$ (without the constant polynomial), the only
such sequences satisfying these conditions are the \emph{exceptional} $X_{1}%
$-Laguerre and $X_{1}$-Jacobi polynomials. Subsequently, during the past five
years, several mathematicians and physicists have discovered and studied other
exceptional orthogonal polynomials $\{p_{n}\}_{n\in\mathbb{N}_{0}\diagdown
A},$ where $A$ is a finite subset of the non-negative integers $\mathbb{N}%
_{0}$ and where $\deg p_{n}=n$ for all $n\in\mathbb{N}_{0}\diagdown A.$ We
call such a sequence an exceptional $X_{\left\vert A\right\vert }$ sequence,
where $\left\vert A\right\vert $ denotes the cardinality of $A.$ All
exceptional sequences found, to date, have the remarkable feature that they
form a complete orthogonal set in their natural Hilbert space setting.

Among the exceptional sets already known are two types of exceptional Laguerre
polynomials, called the Type I and Type II $X_{m}$-Laguerre polynomials, each
omitting $m$ polynomials. In this paper, we briefly discuss these polynomials
and construct the self-adjoint operators generated by their corresponding
second-order differential expressions in the appropriate Hilbert spaces. In
addition, we present a new Type III family of $X_{m}$-Laguerre polynomials
along with a detailed disquisition of its properties. We include several
representations of these polynomials, orthogonality, norms, completeness, the
location of their local extrema and roots, root asymptotics, as well as a
complete spectral study of the second-order Type III exceptional $X_{m}%
$-Laguerre differential expression.

\end{abstract}
\maketitle

\section{Introduction}

An \textit{exceptional orthogonal polynomial system} is a sequence
$\{p_{n}\}_{n\in\mathbb{N}_{0}\diagdown A}$ with the following characteristic properties:

\begin{enumerate}
\item[(a)] $\deg(p_{n})=n$ for $n\in\mathbb{N}_{0}\backslash A,$ where $A$ is
a finite subset of $\mathbb{N}_{0};$

\item[(b)] there exists an interval $I=(a,b)$ and a Lebesgue measurable weight
$w>0$ on $I$ such that%
\[
\int_{I}p_{n}p_{m}w=k_{n}\delta_{n,m}\quad(n,m\in\mathbb{N}_{0}\backslash A)
\]
for some $k_{n}>0$; here $\delta_{n,m}$ denotes the Kronecker delta symbol;

\item[(c)] there exists a second-order differential expression
\[
\ell\lbrack y](x)=a_{2}(x)y^{\prime\prime}(x)+a_{1}(x)y^{\prime}%
(x)+a_{0}(x)y(x)
\]
and, for each $n\in\mathbb{N}_{0}\backslash A,$ there exists a $\lambda_{n}%
\in\mathbb{C}$ such that $y=p_{n}(x)$ is a solution of
\[
\ell\lbrack y](x)=\lambda_{n}y(x)\quad(x\in I);
\]

\item[(d)] for $n\in A,$ there does \underline{not} exist a polynomial $p$ of
degree $n$ such that $y=p(x)$ satisfies $\ell\lbrack y]=\lambda y$ for any
choice of $\lambda\in\mathbb{C};$

\item[(e)] all of the moments%
\[
\int_{I}x^{n}w(x)dx\quad(n\in\mathbb{N}_{0})
\]
of $w$ exist and are finite.
\end{enumerate}

If $\left\vert A\right\vert $ denotes the cardinality of the set $A,$ we call
a sequence $\{p_{n}\}_{n\in\mathbb{N}_{0}\backslash A}$ satisfying conditions
(a)-(e) above an \textit{exceptional sequence of codimension }$\left\vert
A\right\vert $ or an $X_{\left\vert A\right\vert }$\textit{ polynomial
sequence}$.$ The case $A=\{0\}$ was treated in \cite{KMUG}, where the authors
classified $X_{1}$ exceptional orthogonal polynomials and introduced the
$X_{1}$-Laguerre and the $X_{1}$-Jacobi polynomials, so named because of their
similarity to their classical cousins. The fact that these sequences omit a
constant polynomial distinguishes their characterization from the Bochner
classification \cite{Bochner} characterizing the Jacobi, Laguerre, and Hermite
polynomials; of course, the Bochner classification corresponds to
$A=\varnothing.$ Since 2009, several authors have generalized the results of
Kamran, Milson and G\'{o}mez-Ullate in \cite{KMUG} by finding other sequences
of exceptional polynomials $\{p_{n}\}_{n\in\mathbb{N}_{0}\diagdown A},$ where
$A$ is a finite subset of $\mathbb{N}_{0},$ satisfying each of the conditions
in (a)-(e).

In this paper, we discuss three families of exceptional Laguerre polynomials,
each spanning a flag of codimension $m$. Specifically, we deal with two such
exceptional Laguerre sequences associated with
\begin{equation}
A=\{0,1,\ldots,m-1\} \label{A - choice 1}%
\end{equation}
and a new exceptional Laguerre sequence when
\begin{equation}
A=\{1,2,\ldots m\}. \label{A - choice 2}%
\end{equation}
Associated with (\ref{A - choice 1}), the two exceptional Laguerre sequences
are known as the Type I and Type II exceptional Laguerre polynomials,
discovered and introduced into the literature by Odake and Sasaki in
\cite{Odake-Sasaki} and \cite{Odake-Sasaki1}; properties of these two sets
have been studied at length and can be found in, among others, the
contributions \cite{GUKM10, KMUG3, GUMM-Asymptotics}. In this paper, after a
brief review of their properties, we develop the spectral theory of the two
second-order exceptional Laguerre differential equations having the Type I and
Type II sequences as eigenfunctions. As mentioned above, we also present a new
sequence of exceptional Laguerre polynomials, naturally named the Type III
exceptional $X_{m}$-Laguerre polynomials; these polynomials are associated
with the set $A$ given in (\ref{A - choice 2}).

The discovery of the exceptional orthogonal polynomials is one of the most
interesting, and intensive studies, over the past five years in the area of
exactly solvable models in quantum mechanics. We refer the interested reader
to the contributions \cite{KMUG1, GUKM10, KMUG3, KMUG4, KMUG5, KMUG6,
HoOdakeSasaki, HoSasakiZeros2012, Odake-Sasaki, Odake-Sasaki1, Quesne,
Quesne2}, each of which has been influential in our study of this developing
subject. Exceptional orthogonal polynomials and their associated exactly
solvable potentials have applications in a wide range of problems in
mathematical physics, including mass-dependent potentials \cite{Midya-Roy},
supersymmetric quantum mechanics \cite{Dutta-Roy, grandati11}, quasi-exact
solvability \cite{Tanaka}, and Fokker-Planck and Dirac equations \cite{Ho}.
Both mathematicians and physicists have discovered some of these new
orthogonal polynomials as a result of deforming the radial oscillator
potential and the Darboux-P\"{o}schl-Teller potential in terms of an
eigenfunction which is a polynomial of degree $m$. The lowest ($m=1$) such
examples, the exceptional $X_{1}$-Laguerre and $X_{1}$-Jacobi polynomials, are
equivalent to those introduced by Kamran, G\'{o}mez-Ullate and Milson in their
seminal paper \cite{KMUG}. Subsequently, their results were reformulated in
the framework of quantum mechanics and shape-invariant potentials by Quesne
(see \cite{Quesne, Quesne2}) who first related these orthogonal polynomials to
the Darboux transformation and classical orthogonal polynomials. The paper
\cite{Odake-Sasaki} first introduced higher order codimension exceptional
orthogonal polynomials for arbitrary positive integers $m.$ Besides the
Darboux transformation being intimately connected with the derivation of such
exceptional families so is the notion of bispectrality and other tools that
appear in the theory of integrable systems. In mathematical physics, these
functions allow us to write exact solutions to rational extensions of
classical quantum potentials. From the point of view of special functions and
orthogonal polynomials, these exceptional orthogonal polynomials are
polynomial systems formed by solutions to Fuchsian linear equations that
belong to the Heine-Stieltjes class.

From a mathematical point of view, there are several surprises above and
beyond the implications of generalizing the Bochner classification theorem.
Indeed, these new exceptional orthogonal polynomials are wonderful new
examples to illustrate the classical Glazman-Krein-Naimark theory
\cite{Naimark} of constructing self-adjoint operators from Lagrangian
symmetric second-order differential expressions. Furthermore, the fact that
each of the exceptional orthogonal polynomial sequences found, to date, are
complete in their natural Hilbert space setting is nothing short of
remarkable. In particular, it is remarkable that the three exceptional
Laguerre sequences that we discuss in this paper are complete since each of
them are missing $m$ polynomials. On the other hand, their completeness
suggests that interesting M\"{u}ntz-type theorems (see, for example,
\cite{Borwein-Erdelyi} and \cite[Theorem 7.6]{Deutsch}) in weighted $L^{2}%
$-spaces lie in waiting to be discovered.

The contents of this paper are as follows. Section
\ref{Classical Laguerre and Rational Factorizations} reviews some essential
facts about the classical\ Laguerre expression and its solutions. We also give
a brief synopsis of rational factorizations and the Darboux transformation as
it relates to the Laguerre case in this section. In Section
\ref{Type I section}, we review the main properties of the Type I exceptional
$X_{m}$-Laguerre polynomials and then develop the spectral theory of their
associated second-order differential expression. In particular, we construct a
self-adjoint operator, generated from the second-order Type I exceptional
$X_{m}$-Laguerre differential expression, which has the Type I exceptional
$X_{m}$-Laguerre polynomials as eigenfunctions (Theorem
\ref{Type I Self-Adjoint Operator}). We also discuss another interesting
self-adjoint operator (Theorem \ref{Self-Adjoint Operator S Theorem}),
generated by the Type I exceptional $X_{m}$-Laguerre differential expression,
which has a complete set of eigenfunctions involving the Type III exceptional
$X_{m}$-Laguerre polynomials. Section \ref{Type II Section} treats the Type II
exceptional $X_{m}$-Laguerre polynomials in a similar fashion. In Section
\ref{Type III Section}, we introduce the Type III exceptional $X_{m}$-Laguerre
polynomials and we develop many of their properties in full detail. Among
these properties, we will establish several explicit representations of these
polynomials. In Section \ref{Type III Section}, we will also compute the norms
of these polynomials (Theorem \ref{Norms of Type III Polynomials}) in the
appropriate Hilbert space $H$ and show that the sequence of Type III
exceptional $X_{m}$-Laguerre polynomials, despite missing polynomials $p$ of
degrees $1\leq\mathrm{\deg}(p)\leq m$, forms a complete orthogonal set of
polynomials (Theorem \ref{Completeness}) in $H.$ We also determine the
location of the roots of these polynomials (Theorem \ref{Theorem Interlacing}%
); more precisely, we show that the Type III exceptional $X_{m}$-Laguerre
polynomial of degree $m+k$ has $k$ positive roots and $m$ negative roots and
we give careful interlacing properties of these roots with the roots of the
two classical Laguerre polynomials $L_{k-1}^{\alpha+1}(x)$ and $L_{m}%
^{-\alpha-1}(-x)$. We also discuss the asymptotic behavior of the roots of
these Type III exceptional $X_{m}$-Laguerre polynomials (Theorem
\ref{Theorem Asymptotic}). In Section \ref{Type III Section}, we develop
spectral properties of the second-order Type III exceptional $X_{m}$-Laguerre
differential expression $\ell_{m}^{\text{III,}\alpha}[\cdot]$ and, in
particular, determine the self-adjoint operator in $H,$ generated by $\ell
_{m}^{\text{III,}\alpha}[\cdot],$ which has the Type III exceptional $X_{m}%
$-Laguerre polynomials as eigenfunctions (Theorem
\ref{Type III Self-Adjoint Operator}). Lastly, in the Appendix, we list some
examples of Type III exceptional $X_{m}$-Laguerre polynomials.\medskip

\underline{Notation}: Let $\mathbb{N}$ denote the set of positive integers and
$\mathbb{N}_{0}:=\mathbb{N\;}$ $\mathbb{\cup}$ $\{0\}.$ The set $\mathcal{P}$
will denote the vector space of all complex-valued polynomials $p(x)$ in the
real variable $x.$ For $n\in\mathbb{N}_{0},$ let $\mathcal{P}_{n}$ denote the
$(n+1)$-dimensional vector space of all polynomials of degree $\leq n.$

\section{The Classical Laguerre Differential Expression and Rational
Factorizations\label{Classical Laguerre and Rational Factorizations}}

Let
\begin{equation}
\ell^{\alpha}[y]=-xy^{\prime\prime}+(-\alpha-1+x)y^{\prime}
\label{Classical Laguerre}%
\end{equation}
denote the classical Laguerre differential expression. Of course, for each
$m\in\mathbb{N}_{0},$ $y=L_{m}^{\alpha}(x)$ is a solution of
\[
\ell^{\alpha}[y](x)=my(x).
\]

\begin{remark}
In the contributions \cite{KMUG, KMUG1, GUKM10, KMUG3}, the authors define the
Laguerre expression as
\[
\widetilde{\ell}^{\alpha}[y]=xy^{\prime\prime}+(\alpha+1-x)y^{\prime};
\]
for operator-theoretic and spectral-analytic reasons, we elect to define the
Laguerre expression as in \eqref{Classical Laguerre}.
\end{remark}

A rational factorization of $-\ell^{\alpha}[\cdot]$ is an identity of the form%
\begin{equation}
-\ell^{\alpha}=BA-\lambda, \label{Rational Factorization}%
\end{equation}
where $A$ and $B$ are first-order linear differential expressions with
rational coefficients. We call%
\begin{equation}
\widehat{\ell}^{\alpha}:=AB-\lambda\label{Darboux Transformation}%
\end{equation}
the \textit{partner operator} corresponding to the above rational
factorization. Suppose $\phi(x)$ is a quasi-rational solution (that is,
$\dfrac{\phi^{\prime}(x)}{\phi(x)}$ is a rational function) of $A[y]=0.$
Notice, from (\ref{Rational Factorization}), that%
\[
\ell^{\alpha}[\phi](x)=\lambda\phi(x).
\]
The operators $A$ and $B$ are given by%
\[
A[y](x)=b(x)\left(  y^{\prime}(x)-\frac{\phi^{\prime}(x)}{\phi(x)}y(x)\right)
\]
and%
\[
B[y](x)=\widehat{b}(x)\left(  y^{\prime}(x)-\widehat{w}(x)y(x)\right)  ,
\]
where $b(x)$ is a rational function, called the \textit{factorization gauge},
and%
\[
\widehat{b}(x)=\frac{x}{b(x)}\text{ },\text{ }\widehat{w}(x)=-\frac
{\phi^{\prime}(x)}{\phi(x)}+\frac{b^{\prime}(x)}{b(x)}-\frac{1+\alpha-x}{x}.
\]
In this case, the second-order partner operator is given by%
\[
\widehat{\ell}^{\alpha}[y](x):=xy^{\prime\prime}(x)+\widehat{q}(x)y^{\prime
}(x)+\widehat{r}(x)y(x),
\]
where%
\begin{align*}
\widehat{q}(x)  &  =2+\alpha-x-2x\frac{b^{\prime}(x)}{b(x)}\medskip\\
\widehat{r}(x)  &  =-x\left(  \left(  \widehat{w}(x)\right)  ^{\prime}+\left(
\widehat{w}(x)\right)  ^{2}\right)  -\widehat{q}(x)\widehat{w}(x)+\lambda.
\end{align*}
In choosing the factorization gauge $b(x),$ we are guided by two principles:
(a) we want polynomial eigenfunctions of the partner operator and (b) we do
not want these polynomial eigenfunctions to have a common factor. For further 
information on rational factorizations and the above
transformation formulas, we refer the reader to \cite[Section 3]{GUKM10}.

As discussed in \cite{bateman}, the quasi-rational solutions of the
differential equation $\ell^{\alpha}[y]=\lambda y$ are given by the following
four families:
\begin{align}
\phi_{0}(x)  &  =L_{m}^{\alpha}(x) & \lambda &  =m\nonumber\\
\phi_{1}(x)  &  =e^{x}L_{m}^{\alpha}(-x) & \lambda &  =-(\alpha
+1+m)\label{eq:seed4}\\
\phi_{2}(x)  &  =x^{-\alpha}L_{m}^{-\alpha}(x) & \lambda &  =m-\alpha
\label{eq:seed5}\\
\phi_{3}(x)  &  =x^{-\alpha}e^{x}L_{m}^{-\alpha}(-x) & \lambda &  =-(m+1).
\label{eq:seed6}%
\end{align}
Our choice of labels differs from \cite{bateman} to better conform to the Type
$\text{I}$, $\text{II}$, $\text{III}$ nomenclature for the exceptional
Laguerre polynomials.

As described in \cite{GUKM10}, each of these quasi-rational solutions
corresponds to a rational factorization of $\ell^{\alpha}[\cdot]$ and, through
the Darboux transform (\ref{Darboux Transformation}), lead to the Type I, Type
II and, as we will see in this manuscript, the Type III exceptional $X_{m}%
$-Laguerre operators.

\section{The Type I Exceptional $X_{m}$-Laguerre
Polynomials\label{Type I section}}

Unless otherwise indicated, for the Type I exceptional $X_{m}$-Laguerre
polynomials, we assume that the parameter $\alpha>0$. Explicit details of
properties of these polynomials can be found in \cite{GUKM10, KMUG3,
GUMM-Asymptotics}. In this section, we will briefly discuss some properties of
these polynomials and then study the spectral properties of the associated
second-order Type I differential expression.

\subsection{Properties of the Type I Exceptional $X_{m}$-Laguerre Polynomials}

For fixed $m\in\mathbb{N},$ the classical Laguerre polynomial $L_{m}%
^{\alpha-1}(-x)$ has no roots in $[0,\infty)$ when $\alpha>0.$ Taking
$\phi_{1}(x)$ in (\ref{eq:seed4}) as the quasi-rational solution and the
classical Laguerre polynomial $L_{m}^{\alpha}(-x)$ as the factorization gauge,
it can be seen that the classical Laguerre expression $\ell^{\alpha}[\cdot],$
given in (\ref{Classical Laguerre}), may be rewritten as%
\[%
\begin{array}
[c]{l}%
\quad\quad\;-\ell^{\alpha}\;\,=B_{m}^{\text{I,}\alpha}\circ A_{m}%
^{\text{I,}\alpha}+\alpha+m+1,\quad\text{where\medskip}\\
A_{m}^{\text{I,}\alpha}[y](x):=L_{m}^{\alpha}(-x)y^{\prime}(x)-L_{m}%
^{\alpha+1}(-x)y(x)\quad\text{and\medskip}\\
B_{m}^{\text{I,}\alpha}[y](x):=\frac{xy^{\prime}(x)+(1+\alpha)y(x)}%
{L_{m}^{\alpha}(-x)}\,.
\end{array}
\]

\noindent With these definitions, the Type I exceptional $X_{m}$-Laguerre
expression $\ell_{m}^{\text{I,}\alpha}[\cdot]$ may be written as
\[
\ell_{m}^{\text{I,}\alpha}=-\left(  A_{m}^{\text{I},\alpha-1}\circ
B_{m}^{\text{I},\alpha-1}+\alpha+m\,\right)  .
\]
Written out, this expression is given by
\begin{align}
\ell_{m}^{\text{I,}\alpha}[y]:  &  =-\ell^{\alpha}[y](x)+2(\log L_{m}%
^{\alpha-1}(-x))^{\prime}(xy^{\prime}(x)+\alpha
y(x))-my(x)\label{Type I Expression}\\
&  =-xy^{\prime\prime}(x)+\left(  x-\alpha-1+2x\frac{(L_{m}^{\alpha
-1}(-x))^{\prime}}{L_{m}^{\alpha-1}(-x)}\right)  y^{\prime}(x)\nonumber\\
&  \text{\qquad}+\left(  \frac{2\alpha(L_{m}^{\alpha-1}(-x))^{\prime}}%
{L_{m}^{\alpha-1}(-x)}-m\right)  y(x).\nonumber
\end{align}

The Type I exceptional $X_{m}$-Laguerre polynomial $y=L_{m,n}^{\text{I}%
,\alpha}(x)$ $(n\in\mathbb{N}\setminus\{0,1,2,\ldots,m-1\})$ satisfies the
second-order differential equation
\[
\ell_{m}^{\text{I,}\alpha}[y]=\lambda_{n}y\quad(0<x<\infty),
\]
where
\[
\lambda_{n}:=n-m\quad(n=m,m+1,m+2,\ldots)\,
\]
so $\{\lambda_{n}\mid n=m,m+1,m+2,\ldots\}=\mathbb{N}_{0}.$

The $n^{th}$ degree Type I exceptional $X_{m}$-Laguerre polynomial can be
expressed in terms of classical Laguerre polynomials through the formula
\begin{equation}
L_{m,n}^{\text{I},\alpha}(x)=L_{m}^{\alpha}(-x)L_{n-m}^{\alpha-1}%
(x)+L_{m}^{\alpha-1}(-x)L_{n-m-1}^{\alpha}(x)\quad(n\geq m).
\label{Type I Polynomials}%
\end{equation}
In fact, (\ref{Type I Polynomials}) follows by expanding the expression
\[
L_{m,n}^{\text{I},\alpha}(x):=-A_{m}^{\text{I},\alpha-1}[L_{n-m}^{\alpha
-1}](x)\,\quad(n\geq m).
\]
The Type I exceptional $X_{m}$-Laguerre polynomials $\{L_{m,n}^{\text{I}%
,\alpha}\}_{n=m}^{\infty}$ are orthogonal on $(0,\infty)$ with respect to the
weight function
\[
W_{m}^{\text{I},\alpha}(x)=\frac{x^{\alpha}e^{-x}}{(L_{m}^{\alpha-1}(-x))^{2}%
}\qquad(0<x<\infty).
\]

\begin{remark}
\label{Bound of Type I weight}Notice that, since $L_{m}^{\alpha-1}(-x)$ is
positive and increasing on $(0,\infty),$ the function $1/(L_{m}^{\alpha
-1}(-x))^{2}$ is both bounded and bounded away from zero on $(0,\infty).$
Consequently, all moments of the weight function $W_{m}^{\textnormal{I},\alpha}$ on
the interval $(0,\infty)$ exist and are finite.
\end{remark}

In \cite{KMUG3}, the Type I exceptional $X_{m}$-polynomials $\left\{
L_{m,n}^{\text{I,}\alpha}\right\}  _{n=m}^{\infty}$ are shown to be complete
in the Hilbert space $L^{2}((0,\infty);W_{m}^{\text{I},\alpha}).$ With
$\left\Vert \cdot\right\Vert _{m}^{\text{I,}\alpha}$ denoting the norm%
\begin{equation}
\left\Vert f\right\Vert _{m}^{\text{I},\alpha}=\left(  \int_{0}^{\infty
}\left\vert f(x)\right\vert ^{2}W_{m}^{\text{I},\alpha}(x)dx\right)
^{1/2}\quad(f\in L^{2}((0,\infty);W_{m}^{\text{I},\alpha}))
\label{norm for Type I Laguerre}%
\end{equation}
in $L^{2}((0,\infty);W_{m}^{\text{I},\alpha}),$ derived from the inner product%
\[
(f,g)_{m}^{\text{I},\alpha}:=\int_{0}^{\infty}f(x)\overline{g}(x)W_{m}%
^{\text{I},\alpha}(x)dx\quad(f,g\in L^{2}((0,\infty);W_{m}^{\text{I},\alpha
})),
\]
the explicit norms of the Type I exceptional $X_{m}$-Laguerre polynomials are
given by%
\[
\left(  \left\Vert L_{m,n}^{\text{I},\alpha}\right\Vert _{m}^{\text{I},\alpha
}\right)  ^{2}=\frac{(\alpha+n)\Gamma(\alpha+n-m)}{(n-m)!}\quad(n\geq m);
\]
see \cite{GUKM10}.

In \cite{GUMM-Asymptotics}, the authors prove the following two theorems
concerning the zeros of $\{L_{m,n}^{\text{I},\alpha}\}_{n=m}^{\infty}.$

\begin{theorem}
\cite[Proposition 3.2]{GUMM-Asymptotics} For $\alpha>0,$ the Type \text{I}
exceptional $X_{m}$-Laguerre polynomial $L_{m,m+k}^{\textnormal{I},\alpha}(x)$ has
$k$ simple zeros in $(0,\infty)$ and $m$ simple zeros in $(-\infty,0).$ More
specifically, the positive roots of $L_{m,m+k}^{\textnormal{I},\alpha}(x)$ are
located between consecutive roots of $L_{k}^{\alpha}(x)$ and $L_{k-1}^{\alpha
}(x)$ with the smallest positive root of $L_{m,m+k}^{\textnormal{I},\alpha}(x)$
located to the left of the smallest root of $L_{k}^{\alpha}(x).$ The negative
roots of $L_{m,m+k}^{\textnormal{I},\alpha}(x)$ are located between the
consecutive roots of $L_{m-1}^{\alpha}(-x)$ and $L_{m}^{\alpha}(-x).$
\end{theorem}

\begin{theorem}
\cite[Corollary 3.1 and Proposition 3.4]{GUMM-Asymptotics} For $k\geq1,$ the
following asymptotic results for the roots of $L_{m,m+k}^{\textnormal{I},\alpha}(x)$ hold:

\begin{enumerate}
\item[$(a)$] Let $\{j_{\alpha,i}\}$ denote the sequence of zeros of the Bessel
function of the first kind $J_{\alpha}(x)$ listed in increasing order and let
$\{x_{k,i}^{\alpha}\}_{i=1}^{k}$ denote the positive zeros of $L_{m,m+k}%
^{\textnormal{I},\alpha}(x)$ arranged in increasing order$.$ Then%
\[
\lim_{k\rightarrow\infty}kx_{k,i}^{\alpha}=\frac{j_{\alpha,i}^{2}}{4}%
\quad(i\in\mathbb{N}).
\]

\item[$(b)$] As $k\rightarrow\infty,$ the $m$ negative roots of $L_{m,m+k}%
^{\textnormal{I},\alpha}(x)$ converge to the $m$ roots of $L_{m}^{\alpha-1}(-x).$
\end{enumerate}
\end{theorem}

\subsection{Type I Exceptional $X_{m}$-Laguerre Spectral
Analysis\label{Type I Spectral Analysis}}

In Lagrangian symmetric form, the Type I exceptional $X_{m}$-Laguerre
differential expression (\ref{Type I Expression}) is given by
\begin{align}
\ell_{m}^{\text{I,}\alpha}[y](x)=\frac{1}{W_{m}^{\text{I},\alpha}(x)}  &
\left(  -\left(  \frac{x^{\alpha+1}e^{-x}}{(L_{m}^{\alpha-1}(-x))^{2}%
}y^{\prime}(x)\right)  ^{\prime}\right. \label{Type I Symmetric Form}\\
&  +\left.  \left(  \frac{2\alpha x^{\alpha}e^{-x}(L_{m}^{\alpha
-1}(-x))^{\prime}}{(L_{m}^{\alpha-1}(-x))^{3}}-\frac{mx^{\alpha}e^{-x}}%
{(L_{m}^{\alpha-1}(-x))^{2}}\right)  y(x)\right)  .\nonumber
\end{align}
When $m=1,$ the spectral analysis of (\ref{Type I Symmetric Form}) in
$L^{2}((0,\infty);W_{1}^{\text{I},\alpha})$ was completed in
\cite{Atia-Littlejohn-Stewart}.

The maximal domain associated with $\ell_{m}^{\text{I},\alpha}[\cdot]$ in the
Hilbert space $L^{2}\left(  (0,\infty);W_{m}^{\text{I},\alpha}\right)  $ is
defined to be
\[
\Delta_{m}^{\text{I,}\alpha}:=\left\{  f:(0,\infty)\rightarrow\mathbb{C}\mid
f,f^{\prime}\in AC_{\text{\textrm{loc}}}(0,\infty);f,\ell_{m}^{\text{I,}%
\alpha}[f]\in L^{2}((0,\infty);W_{m}^{\text{I},\alpha})\right\}  .
\]
The associated maximal operator
\[
T_{1,m}^{\text{I},\alpha}\ :\mathcal{D}(T_{1,m}^{\text{I,}\alpha})\subset
L^{2}((0,\infty);W_{m}^{\text{I},\alpha})\rightarrow L^{2}((0,\infty
);W_{m}^{\text{I},\alpha}),
\]
is defined to be
\begin{align*}
T_{1,m}^{\text{I,}\alpha}\ f  &  =\ell_{m}^{\text{I,}\alpha}[f]\\
f\in\mathcal{D}(T_{1,m}^{\text{I,}\alpha}):  &  =\Delta_{m}^{\text{I,}\alpha}.
\end{align*}

For $f,g\in\Delta_{m}^{\text{I,}\alpha}$, Green's formula is given by
\[
\int_{0}^{\infty}\ell_{m}^{\text{I,}\alpha}[f](x)\overline{g}(x)W_{m}%
^{\text{I},\alpha}(x)dx=[f,g]_{m}^{\text{I,}\alpha}(x)\mid_{x=0}^{x=\infty
}+\int_{0}^{\infty}f(x)\ell_{m}^{\text{I,}\alpha}[\overline{g}](x)W_{m}%
^{\text{I},\alpha}(x)dx\,,
\]
where $[\cdot,\cdot]_{m}^{\text{I,}\alpha}(\cdot)$ is the sesquilinear form
defined by
\[
\left[  f,g\right]  _{m}^{\text{I},\alpha}(x):=\frac{x^{\alpha+1}e^{-x}%
}{(L_{m}^{\alpha-1}(-x))^{2}}(f(x)\bar{g}^{\prime}(x)-f^{\prime}(x)\bar
{g}(x))\quad(0<x<\infty),
\]
and where
\[
\left[  f,g\right]  _{m}^{\text{I},\alpha}(x)\mid_{x=0}^{x=\infty}:=\left[
f,g\right]  _{m}^{\text{I},\alpha}(\infty)-\left[  f,g\right]  _{m}%
^{\text{I},\alpha}(0)\,.
\]
By Green's formula and the definition of $\Delta_{m}^{\text{I,}\alpha},$ both
limits
\[
\left[  f,g\right]  _{m}^{\text{I},\alpha}(\infty):=\lim_{x\rightarrow\infty
}[f,g]_{m}^{\text{I,}\alpha}(x)\text{ and }\left[  f,g\right]  _{m}%
^{\text{I},\alpha}(0):=\lim_{x\rightarrow0^{+}}[f,g]_{m}^{\text{I,}\alpha}(x)
\]
exist and are finite for all $f,g\in\Delta_{m}^{\text{I,}\alpha}.$

The adjoint of the maximal operator in $L^{2}\left(  (0,\infty);W_{m}%
^{\text{I},\alpha}\right)  $ is the minimal operator
\[
T_{0,m}^{\text{I},\alpha}\ :\mathcal{D}(T_{0,m}^{\text{I},\alpha})\subset
L^{2}((0,\infty);W_{m}^{\text{I},\alpha})\rightarrow L^{2}((0,\infty
);W_{m}^{\text{I},\alpha}),
\]
defined by%
\begin{align*}
T_{0,m}^{\text{I},\alpha}f  &  =\ell_{m}^{\text{I,}\alpha}[f]\\
f\in\mathcal{D}(T_{0,m}^{\text{I},\alpha}):  &  =\{f\in\Delta_{m}%
^{\text{I,}\alpha}\mid\lbrack f,g]_{m}^{\text{I,}\alpha}\mid_{x=0}^{x=\infty
}=0\text{ for all }g\in\Delta_{m}^{\text{I,}\alpha}\}.
\end{align*}

We seek to find the self-adjoint extension $T_{m}^{\text{I},\alpha}\ $in
$L^{2}\left(  (0,\infty);W_{m}^{\text{I},\alpha}\right)  ,$ generated by
$\ell_{m}^{\text{I,}\alpha}[\cdot]$, which has the Type I exceptional $X_{m}%
$-Laguerre polynomials $\left\{  L_{m,n}^{\text{I},\alpha}\right\}
_{n=m}^{\infty}$ as eigenfunctions. To do this, we first need to study the
behavior of solutions near the singular endpoints $x=0$ and $x=\infty$ in
order to determine the deficiency indices and to determine the appropriate
boundary conditions (if any).

Our first goal is to show the following theorem:

\begin{theorem}
\label{LP/LC Type I Theorem}For $\alpha>0$, let $\ell_{m}^{\textnormal{I,}\alpha
}[\cdot]$ be the Type I exceptional $X_{m}$-Laguerre differential expression
\eqref{Type I Expression} on the interval $(0,\infty)$.

\begin{enumerate}
\item[$(a)$] $\ell_{m}^{\textnormal{I,}\alpha}[\cdot]$ is in the limit-circle case
at $x=0$ when $0<\alpha<1$ and is in the limit-point case at $x=0$ when
$\alpha\geq1$.

\item[$(b)$] $\ell_{m}^{\textnormal{I,}\alpha}[\cdot]$ is in the limit-point case
at $x=\infty$ for any choice of $\alpha>0$.
\end{enumerate}

\begin{proof}
$(a)$: The endpoint $x=0$ is, in the sense of Frobenius, a regular singular
endpoint of the Type I exceptional $X_{m}$-Laguerre expression $\ell
_{m}^{\text{I,}\alpha}[y]=0$. The Frobenius indicial equation at $x=0$ is
\[
r(r+\alpha)=0\,.
\]
Consequently, two linearly independent solutions of $\ell_{m}^{\text{I,}%
\alpha}[y]=0$ on $(0,\infty)$ will behave asymptotically like
\[
z_{1}(x)=1\quad\text{and}\quad z_{2}(x)=x^{-\alpha}%
\]
near $x=0$. Now, for any $\alpha>0$, we see from Remark
\ref{Bound of Type I weight} that
\[
\int_{0}^{\infty}\left\vert z_{1}(x)\right\vert ^{2}W_{m}^{\text{I},\alpha
}(x)\,dx=\int_{0}^{\infty}\frac{x^{\alpha}e^{-x}}{(L_{m}^{\alpha-1}(-x))^{2}%
}\,dx<\infty\,.
\]
However, for any choice of $x^{\ast}\in(0,\infty)$,
\[
\int_{0}^{x^{\ast}}\left\vert z_{2}(x)\right\vert ^{2}W_{m}^{\text{I},\alpha
}(x)\,dx=\int_{0}^{x^{\ast}}\frac{x^{-\alpha}e^{-x}}{(L_{m}^{\alpha
-1}(-x))^{2}}\,dx<\infty
\]
only when $0<\alpha<1$. In the language of the Weyl limit-point/limit-circle
theory, it follows that the Type I exceptional $X_{m}$-Laguerre differential
expression is in the limit-circle case at $x=0$ when $0<\alpha<1$ and is in
the limit-point case at $x=0$ when $\alpha\geq1$.\medskip

$(b)$: Since $x=\infty$ is an irregular singular endpoint of the Type I
exceptional $X_{m}$-Laguerre differential expression, the above Frobenius
method cannot be employed. Fortunately, we are able to explicitly solve the
differential equation
\[
\ell_{m}^{\text{I,}\alpha}[y](x)=0\quad(0<x<\infty)
\]
for a basis $\left\{  y_{1}(x),y_{2}(x)\right\}  $ of solutions and, from
this, we are able to determine the $L^{2}$ behavior of these solutions near
$x=\infty.$ The function $y_{1}(x)=L_{m,m}^{\text{I},\alpha}(x)=L_{m}^{\alpha
}(-x)$, the Type I exceptional $X_{m}$-Laguerre polynomial of degree $m$, is
one solution of $\ell_{m}^{\text{I},\alpha}[y](x)=0$ on $(0,\infty).$ Using
the well-known reduction of order method, we obtain a second linearly
independent solution $y_{2}(x)$. Indeed,
\begin{equation}
y_{2}(x)=L_{m}^{\alpha}(-x)\int_{a}^{x}\frac{e^{t}(L_{m}^{\alpha-1}(-t))^{2}%
}{t^{\alpha+1}(L_{m}^{\alpha}(-t))^{2}}\,dt\,, \label{y_2}%
\end{equation}
where $a$ is a fixed, arbitrary, positive constant. Clearly $y_{1}$ $\in
L^{2}((0,\infty);W_{m}^{\text{I},\alpha})$. However, as we now show,
$y_{2}\notin L^{2}((0,\infty);W_{m}^{\text{I},\alpha})$. To see this, note
that since
\[
\lim_{t\rightarrow\infty}\left(  \frac{L_{m}^{\alpha-1}(-t)}{L_{m}^{\alpha
}(-t)}\right)  ^{2}=1\,,
\]
there exists $x_{0}>0$ such that%
\begin{equation}
\left(  \frac{L_{m}^{\alpha-1}(-t)}{L_{m}^{\alpha}(-t)}\right)  ^{2}\geq
\frac{1}{2}\quad(t\geq x_{0}).
\label{Estimate of square of ratios of Laguerre}%
\end{equation}
Hence, for large enough $x_{1}\geq x_{0}$
\begin{align*}
\int_{x_{0}}^{x}\frac{e^{t}}{t^{\alpha+1}}\left(  \frac{L_{m}^{\alpha-1}%
(-t)}{L_{m}^{\alpha}(-t)}\right)  ^{2}\,dt  &  \geq\frac{1}{2}\int_{x_{0}}%
^{x}\frac{e^{t}}{t^{\alpha+1}}\,dt\\
&  \geq\frac{A}{2}\int_{x_{1}}^{x}e^{t/2}\,dt\quad\text{where }A=\frac
{e^{\alpha+1}}{2^{\alpha+1}(\alpha+1)^{\alpha+1}}\quad(x\geq x_{1}).
\end{align*}
Moreover, for large enough $x_{2}\geq x_{1},$
\[
\frac{A}{2}\int_{x_{2}}^{x}e^{t/2}\,dt\geq\frac{A}{2}e^{x/2}\text{ \quad
}(x\geq x_{2}).
\]
Hence, from $($\ref{y_2}$)$ with the choice $a=x_{2},$ we see that
\begin{align}
\left\vert y_{2}(x)\right\vert ^{2}  &  =\left(  L_{m}^{\alpha}(-x)\right)
^{2}\left(  \int_{x_{2}}^{x}\frac{e^{t}}{t^{\alpha+1}}\left(  \frac
{L_{m}^{\alpha-1}(-t)}{L_{m}^{\alpha}(-t)}\right)  ^{2}\,dt\right)
^{2}\label{y_2 estimate}\\
&  \geq\left(  \frac{A}{2}\right)  ^{2}\left(  L_{m}^{\alpha}(-x)\right)
^{2}e^{x}\quad(x\geq x_{2}).\nonumber
\end{align}
Mimicking the argument in $($\ref{Estimate of square of ratios of Laguerre}%
$),$ there exists an $x_{3}\geq x_{2}$ such that
\[
\left(  \frac{L_{m}^{\alpha}(-t)}{L_{m}^{\alpha-1}(-t)}\right)  ^{2}\geq
\frac{1}{2}\quad(t\geq x_{3}).
\]
Consequently, from $($\ref{y_2 estimate}$),$ we see that for $t\geq x_{3},$
\begin{align*}
&  \int_{x_{3}}^{\infty}\left\vert y_{2}(t)\right\vert ^{2}W_{m}%
^{\text{I},\alpha}(t)dt\\
&  \geq\left(  \frac{A}{2}\right)  ^{2}\int_{x_{3}}^{\infty}\,t^{\alpha
}\left(  \frac{L_{m}^{\alpha}(-t)}{L_{m}^{\alpha-1}(-t)}\right)  ^{2}dt\\
&  \geq\frac{A^{2}}{8}\int_{x_{3}}^{\infty}t^{\alpha}dt=\infty.
\end{align*}
Therefore, $y_{2}\notin L^{2}((x_{3},\infty);W_{m}^{\text{I},\alpha})$.
\end{proof}
\end{theorem}

As a result,

\begin{theorem}
Let $T_{0,m}^{\textnormal{I},\alpha}\ $be the minimal operator in $L^{2}\left(
(0,\infty);W_{m}^{\textnormal{I},\alpha}\right)  $ generated by the Type I
exceptional $X_{m}$-Laguerre differential expression $\ell_{m}%
^{\textnormal{I,}\alpha}[\cdot]$.

\begin{enumerate}
\item[$(a)$] {If $0<\alpha<1$, the deficiency index of $T_{0,m}%
^{\textnormal{I},\alpha}\ $is $(1,1)$;}

\item[$(b)$] {If $\alpha\geq1$, the deficiency index of $T_{0,m}%
^{\textnormal{I},\alpha}\ $is $(0,0)$.}
\end{enumerate}
\end{theorem}

If $\alpha\geq1,$ there is only one self-adjoint extension (restriction) of
the minimal operator {$T_{0,m}^{\text{I},\alpha}$ (maximal operator
$T_{1,m}^{\text{I,}\alpha}$}), namely $T_{m}^{\text{I},\alpha}:=$
{$T_{0,m}^{\text{I},\alpha}=T_{1,m}^{\text{I,}\alpha}.$ }However, when
$0<\alpha<1$, there are infinitely many self-adjoint extensions of
{$T_{0,m}^{\text{I},\alpha}.$ }Furthermore, when $0<\alpha<1,$ in order to
obtain a self-adjoint extension of the minimal operator $T_{0,m}%
^{\text{I},\alpha}\ $having the Type I exceptional $X_{m}$-Laguerre
polynomials $\left\{  L_{m,n}^{\text{I},\alpha}\right\}  _{n=m}^{\infty}$ as
eigenfunctions, the Glazman-Krein-Naimark theory (see \cite{Naimark}) requires
that we impose one particular boundary condition of the form,
\[
\left[  f,g_{0}\right]  _{m}^{\text{I},\alpha}(0)=0\quad(f\in\Delta
_{m}^{\text{I,}\alpha})\,,
\]
where $g_{0}\in\Delta_{m}^{\text{I,}\alpha}\setminus\mathcal{D}(T_{0,m}%
^{\text{I},\alpha})$. We claim $g_{0}\equiv1$ on $(0,\infty)$ is an
appropriate choice.

Note that the function $y(x)=x^{-\alpha}\in L^{2}\left(  (0,\infty
);W_{m}^{\text{I},\alpha}\right)  $ when $0<\alpha<1$. Remarkably, it is the
case that
\[
\ell_{m}^{\text{I},\alpha}[x^{-\alpha}]=(-m-\alpha)x^{-\alpha}\,;
\]
hence, it follows that $x^{-\alpha}\in\Delta_{m}^{\text{I,}\alpha}$ for
$0<\alpha<1$. Additionally,
\begin{equation}
\lbrack x^{-\alpha},1]_{m}^{\text{I},\alpha}(0)=\alpha\lim_{x\rightarrow0^{+}%
}\frac{e^{-x}}{(L_{m}^{\alpha-1}(-x))^{2}}\neq0\,
\label{No Minimal Domain Function}%
\end{equation}
so $g_{0}=1\in\Delta_{m}^{\text{I,}\alpha}\backslash\mathcal{D}(T_{0,m}%
^{\text{I},\alpha}\ ).$ Furthermore, a calculation shows that
\[
\lbrack f,1]_{m}^{\text{I},\alpha}(0)=0\Longleftrightarrow\lim_{x\rightarrow
0^{+}}x^{\alpha+1}f^{\prime}(x)=0.
\]
Therefore, we obtain the following theorem.

\begin{theorem}
\label{Type I Self-Adjoint Operator}Let $\alpha>0.$

\begin{enumerate}
\item[$(a)$] Suppose $0<\alpha<1$. The operator
\[
T_{m}^{\textnormal{I},\alpha}:\mathcal{D}(T_{m}^{\textnormal{I},\alpha})\subset
L^{2}\left(  (0,\infty);W_{m}^{\textnormal{I},\alpha}\right)  \rightarrow
L^{2}\left(  (0,\infty);W_{m}^{\textnormal{I},\alpha}\right)  \,,
\]
defined by
\begin{align*}
T_{m}^{\textnormal{I},\alpha}\ f  &  =\ell_{m}^{\textnormal{I,}\alpha}[f]\\
f\in\mathcal{D}(T_{m}^{\textnormal{I},\alpha}):  &  =\left\{  f\in\Delta
_{m}^{\textnormal{I,}\alpha}\mid\lim_{x\rightarrow0^{+}}x^{\alpha+1}f^{\prime
}(x)=0\right\}  ,
\end{align*}
is self-adjoint in $L^{2}\left(  (0,\infty);W_{m}^{\textnormal{I},\alpha}\right)
$ and has the Type I exceptional $X_{m}$-Laguerre polynomials $\left\{
L_{m,n}^{\textnormal{I},\alpha}\right\}  _{n=m}^{\infty}$ as eigenfunctions.
Moreover, the spectrum of $T_{m}^{\textnormal{I},\alpha}\ $consists only of
eigenvalues and is given by
\[
\sigma(T_{m}^{\textnormal{I},\alpha}\ )=\mathbb{N}_{0}\,.
\]

\item[$(b)$] Suppose $\alpha\geq1.$ The operator
\[
T_{m}^{\textnormal{I},\alpha}:\mathcal{D}(T_{m}^{\textnormal{I},\alpha})\subset
L^{2}\left(  (0,\infty);W_{m}^{\textnormal{I},\alpha}\right)  \rightarrow
L^{2}\left(  (0,\infty);W_{m}^{\textnormal{I},\alpha}\right)  \,,
\]
defined by
\begin{align*}
T_{m}^{\textnormal{I},\alpha}f  &  =\ell_{m}^{\textnormal{I,}\alpha}[f]\\
f\in\mathcal{D}(T_{m}^{\textnormal{I},\alpha}):  &  =\Delta_{m}^{\textnormal{I,}\alpha
},
\end{align*}
is self-adjoint in $L^{2}\left(  (0,\infty);W_{m}^{\textnormal{I},\alpha}\right)
$ and has the Type I exceptional $X_{m}$-Laguerre polynomials $\left\{
L_{m,n}^{\textnormal{I},\alpha}\right\}  _{n=m}^{\infty}$ as eigenfunctions.
Moreover, the spectrum of $T_{m}^{\textnormal{I},\alpha}\ $consists only of
eigenvalues and is given by
\[
\sigma(T_{m}^{\textnormal{I},\alpha})=\mathbb{N}_{0}\,.
\]

\end{enumerate}
\end{theorem}

\subsection{Spectral Properties of Another Self-Adjoint Operator in
$L^{2}((0,\infty);W_{m}^{\textnormal{I},\alpha})$}

From (\ref{No Minimal Domain Function}) we can interchange the roles of $1$
and $x^{-\alpha}$ to obtain another interesting self-adjoint operator
$S_{m}^{\text{I},\alpha},$ generated by $\ell_{m}^{\text{I,}\alpha}[\cdot],$
in $L^{2}((0,\infty);W_{m}^{\text{I},\alpha})$ with different boundary
conditions. Observing that%
\begin{equation}
0=[f,x^{-\alpha}]_{m}^{\text{I},\alpha}(0)\Longleftrightarrow\lim
_{x\rightarrow0^{+}}(xf^{\prime}(x)+\alpha f(x))=0,
\label{BC for SA operator S}%
\end{equation}
we are now in position to prove the following theorem regarding the
self-adjoint operator $S_{m}^{\text{I},\alpha}$, defined below in
(\ref{SA S Operator}), in the space $L^{2}((0,\infty);W_{m}^{\text{I},\alpha
}).$ This operator $S_{m}^{\text{I},\alpha}$ is only quasi-isospectral to the
classical Laguerre operator in the sense that the commutation transformation
that relates the classical Laguerre operator to its Type III counterpart
represents a state-adding transformation in the sense of Deift \cite{deift78}.

\begin{theorem}
\label{Self-Adjoint Operator S Theorem}Suppose $0<\alpha<1.$ Define
\[
S_{m}^{\textnormal{I},\alpha}:\mathcal{D}(S_{m}^{\textnormal{I},\alpha})\subset
L^{2}((0,\infty);W_{m}^{\textnormal{I},\alpha})\rightarrow L^{2}((0,\infty
);W_{m}^{\textnormal{I},\alpha})
\]
by%
\begin{equation}
\left\{
\begin{array}
[c]{rl}%
S_{m}^{\textnormal{I},\alpha}f & =\ell_{m}^{\textnormal{I,}\alpha}[f]\\
\mathcal{D}(S_{m}^{\textnormal{I},\alpha}) & =\{f\in\Delta_{m}^{\textnormal{I,}\alpha
}\mid\lim_{x\rightarrow0^{+}}(xf^{\prime}(x)+\alpha f(x))=0\}
\end{array}
\right.  \label{SA S Operator}%
\end{equation}
Then $S_{m}^{\textnormal{I},\alpha}$ is self-adjoint in $L^{2}((0,\infty
);W_{m}^{\textnormal{I},\alpha}).$ Furthermore,
\[
\left\{  x^{-\alpha}L_{m,n}^{\textnormal{III},-\alpha}(x)\mid
n=0,m+1,m+2,m+3,\ldots\right\}
\]
forms a complete set of $($orthogonal$)$ eigenfunctions of $S_{m}%
^{\textnormal{I},\alpha}$ in $L^{2}((0,\infty);W_{m}^{\textnormal{I},\alpha}),$ where%
\[
\left\{  L_{m,n}^{\textnormal{III},\alpha}\mid n=0,m+1,m+2,m+3,\ldots\right\}
\quad(-1<\alpha<0)
\]
is the sequence of Type III exceptional $X_{m}$-Laguerre polynomials which we
introduce below in Section \ref{Type III Section}. Finally,%
\[
\sigma(S_{m}^{\textnormal{I},\alpha})=\sigma_{\text{\textrm{p}}}(S_{m}%
^{\textnormal{I},\alpha})=\{n-m-\alpha\mid n=0,m+1,m+2,m+3,\ldots\}.
\]

\begin{proof}
The self-adjointness of $S_{m}^{\text{I},\alpha}$ is clear from the
Glazman-Krein-Naimark theory. We need to prove

\begin{enumerate}
\item[(i)] For $n=0,m+1,m+2,m+3,\ldots,$ $y=x^{-\alpha}L_{m,n}^{\text{III}%
,-\alpha}(x)$ is a solution of
\begin{equation}
\ell_{m}^{\text{I,}\alpha}[y](x)=(n-m-\alpha)y(x).
\label{multiplier function solutions}%
\end{equation}

\item[(ii)] For $n=0,m+1,m+2,m+3,\ldots,$%
\[
x^{-\alpha}L_{m,n}^{\text{III},-\alpha}\in\mathcal{D}(S_{m}^{\text{I},\alpha
}).
\]

\item[(iii)] $\{x^{-\alpha}L_{m,n}^{\text{III},-\alpha}\mid
n=0,m+1,m+2,m+3,\ldots\}$ is a complete $($orthogonal$)$ set of eigenfunctions
of $S_{m}^{\text{I},\alpha}$ in $L^{2}((0,\infty);W_{m}^{\text{I},\alpha});$
equivalently,
\[
\mathrm{span}\left\{  x^{-\alpha}L_{m,n}^{\text{III},-\alpha}\mid
n=0,m+1,m+2,m+3,\ldots\right\}
\]
is dense in $L^{2}((0,\infty);W_{m}^{\text{I},\alpha}).$

\item[(iv)] The spectrum of $S_{m}^{\text{I},\alpha}$ is pure point spectrum
and given by%
\[
\sigma(S_{m}^{\text{I},\alpha})=\{n-m-\alpha\mid n=0,m+1,m+2,m+3,\ldots\}.
\]

\end{enumerate}

Appealing to $($\ref{Type III Laguerre rep}$),$ $($%
\ref{L3 eigenvalue equation}$),$ and $($\ref{L3 eigenvalues}$)$ in Section
\ref{Type III Section}$,$ we see from $($\ref{I to III transformation}$)$ and
$($\ref{I to III transformation 2}$)$ that $y=x^{-\alpha}L_{m,n}%
^{\text{III},-\alpha}(x)$ satisfies
\[
\ell_{m}^{\text{I},\alpha}[y](x)=(n-m-\alpha)y(x)\quad(n=0,m+1,m+2,m+3,\ldots
).
\]
This establishes $($\textrm{i}$).$ From the proof of completeness in part
$($iii$)$ below, part $($iv$)$ will follow.\medskip\ \newline It is clear,
since $0<\alpha<1,$ that $x^{-\alpha}L_{m,n}^{\text{III},-\alpha}\in\Delta
_{m}^{\text{I,}\alpha}.$ Moreover, for fixed $n\in\{0,m+1,m+2,\ldots\},$ let
$h(x)=x^{-\alpha}L_{m,n}^{\text{III},-\alpha}(x);$ a calculation shows that%
\begin{align*}
\lim_{x\rightarrow0^{+}}xh^{\prime}(x)+\alpha h(x)  &  =\lim_{x\rightarrow
0^{+}}\left(  -\alpha x^{-\alpha}L_{m,n}^{\text{III},-\alpha}(x)+x^{-\alpha
+1}\left(  L_{m,n}^{\text{III},-\alpha}(x)\right)  ^{\prime}+\alpha
x^{-\alpha}L_{m,n}^{\text{III},-\alpha}(x)\right) \\
&  =\lim_{x\rightarrow0^{+}}x^{-\alpha+1}\left(  L_{m,n}^{\text{III},-\alpha
}(x)\right)  ^{\prime}=0.
\end{align*}
Hence $h\in\mathcal{D}(S_{m}^{\text{I},\alpha})$ and this proves
$($\textrm{ii}$).\medskip$ \newline To prove $($\textrm{iii}$),$ let $f\in
L^{2}((0,\infty);W_{m}^{\text{I},\alpha})$ and $\varepsilon>0.$ A calculation
shows that%
\begin{equation}
\left\Vert f\right\Vert _{m}^{\text{I},\alpha}=\left\Vert \widetilde{f}%
\right\Vert ^{-\alpha}, \label{SA S 1}%
\end{equation}
where $\widetilde{f}(x)=\dfrac{f(x)}{x^{-\alpha}L_{m}^{\alpha-1}(-x)},$
$\left\Vert \cdot\right\Vert _{m}^{\text{I},\alpha}$ is the norm defined in
$($\ref{norm for Type I Laguerre}$)$ and $\left\Vert \cdot\right\Vert
^{\alpha}$ $(\alpha>-1)$ is the $($usual$)$ norm in $L^{2}((0,\infty
);x^{\alpha}e^{-x})$ given by%
\[
\left\Vert g\right\Vert ^{\alpha}:=\left(  \int_{0}^{\infty}\left\vert
g(x)\right\vert ^{2}x^{\alpha}e^{-x}dx\right)  ^{1/2}\quad(g\in L^{2}%
((0,\infty);x^{\alpha}e^{-x})).
\]
Let%
\[
\eta(x)=L_{m}^{\alpha-1}(-x);
\]
it is clear that $\eta(x)\neq0$ on $(0,\infty)$. Appealing to Lemma
\ref{polynomial x P dense result} in Section \ref{Type III Section}, it
follows that
\[
\left\{  L_{m}^{\alpha-1}(-x)p\mid p\in\mathcal{P}\right\}
\]
is dense in $L^{2}((0,\infty);x^{-\alpha}e^{-x}),$ where $\mathcal{P}$ is the
space of all polynomials in the real variable $x.$ Hence, from $($%
\ref{SA S 1}$),$ there exists $p\in\mathcal{P},$ say \textrm{deg}$(p)=n,$ such
that%
\begin{equation}
\varepsilon>\left\Vert \widetilde{f}-L_{m}^{\alpha-1}(-x)p\right\Vert
^{-\alpha}=\left\Vert f-x^{-\alpha}\left(  L_{m}^{\alpha-1}(-x)\right)
^{2}p\right\Vert _{m}^{\text{I},\alpha}. \label{SA S 3}%
\end{equation}
Define%
\[
\mathcal{E}_{n+2m}:=\text{\textrm{span}}\left\{  x^{-\alpha}L_{m,j}%
^{\text{III},-\alpha}(x)\mid j=0,m+1,m+2,\ldots,n+2m\right\}
\]
and%
\[
\mathcal{F}_{n+2m}:=\left\{  x^{-\alpha}q(x)\mid q\in\mathcal{P}_{n+2m},\text{
}q^{\prime}(-x_{j})=0\text{ }(j=1,2,\ldots,m)\right\}  ,
\]
where $\mathcal{P}_{n+2m}\subset$ $\mathcal{P}$ is the space of all
polynomials of degree $\leq n+2m$ and $\{x_{j}\}_{j=1}^{m}\subset(0,\infty)$
are the simple, distinct roots of the Laguerre polynomial $L_{m}^{\alpha
-1}(x)$. Note that both $\mathcal{E}_{n+2m}$ and $\mathcal{F}_{n+2m}$ are
subspaces of $L^{2}((0,\infty);W_{m}^{\text{I},\alpha})$ with%
\begin{equation}
\text{\textrm{dim}}\left(  \mathcal{F}_{n+2m}\right)  =\text{\textrm{dim}%
}\left(  \mathcal{E}_{n+2m}\right)  =n+m+1. \label{SA S 4}%
\end{equation}
From part $($\textrm{i}$),$ we know that $\mathcal{E}_{n+2m}$ is an invariant
subspace of the Type \textrm{I }exceptional\textrm{ }$X_{m}$-Laguerre
differential expression $\ell_{m}^{\text{I,}\alpha}[\cdot];$ that is to say,%
\[
\ell_{m}^{\text{I,}\alpha}[\mathcal{E}_{n+2m}]\subset\mathcal{E}_{n+2m}.
\]
Furthermore, it is clear that
\begin{equation}
x^{-\alpha}\left(  L_{m}^{\alpha-1}(-x)\right)  ^{2}p\in\mathcal{F}_{n+2m}.
\label{SA S 5}%
\end{equation}
In lieu of $($\ref{SA S 3}$),$ in order to show that
\[
\mathrm{span}\left\{  x^{-\alpha}L_{m,n}^{\text{III},-\alpha}\mid
n=0,m+1,m+2,m+3,\ldots\right\}
\]
is dense in $L^{2}((0,\infty);W_{m}^{\text{I},\alpha}),$ we need to show that
\begin{equation}
x^{-\alpha}\left(  L_{m}^{\alpha-1}(-x)\right)  ^{2}p\in\mathcal{E}_{n+2m}.
\label{SA S 6}%
\end{equation}
\underline{Claim}: $\mathcal{E}_{n+2m}\subset\mathcal{F}_{n+2m}$\newline Let
$g\in$ $\mathcal{E}_{n+2m}$ so $g(x)=x^{-\alpha}q(x)$ where $q(x)$ is a linear
combination of $L_{m,0}^{\text{III},-\alpha}(x),$ $L_{m,m+1}^{\text{III}%
,-\alpha}(x),$ $L_{m,m+2}^{\text{III,}-\alpha}(x),\ldots,$ $L_{m,n}%
^{\text{III},-\alpha}(x).$ A calculation shows%
\begin{align*}
&  \ell_{m}^{\text{I,}\alpha}[g](x)\\
&  =x^{-\alpha}\left(  (\alpha-1)q^{\prime}(x)-xq^{\prime\prime}%
(x)+xq^{\prime}(x)-\alpha q(x)+2x\frac{(L_{m}^{\alpha-1}(-x))^{\prime}}%
{L_{m}^{\alpha-1}(-x)}q^{\prime}(x)-mq(x)\right)  .
\end{align*}
Since $\ell_{m}^{\text{I,}\alpha}[g]\in\mathcal{E}_{n+2m}$ we see that, from
the term
\[
2x\frac{(L_{m}^{\alpha-1}(-x))^{\prime}}{L_{m}^{\alpha-1}(-x)}q^{\prime}(x),
\]
that, necessarily, $q\in\mathcal{P}_{n+2m}$ and $q^{\prime}(-x_{j})=0$ for
$j=1,2,\ldots,m.$ This proves the claim. Notice also that this implies that
the polynomial $q$ is either a constant or \textrm{deg}$(q)\geq m+1.$\newline
By $($\ref{SA S 4}$),$ we see that%
\[
\mathcal{E}_{n+2m}=\mathcal{F}_{n+2m}.
\]
From $($\ref{SA S 5}$),$ we obtain $x^{-\alpha}\left(  L_{m}^{\alpha
-1}(-x)\right)  ^{2}p\in\mathcal{E}_{n+2m},$ establishing $($\ref{SA S 6}$)$
and this proves part $($\textrm{iii}$).$\newline Finally, part $($%
\textrm{iv}$)$ follows immediately from part $($\textrm{iii}$)$ and
$($\ref{multiplier function solutions}$).$ This completes the proof of the theorem.
\end{proof}
\end{theorem}

\section{The Type II Exceptional $X_{m}$-Laguerre
Polynomials\label{Type II Section}}

Throughout this section, we let $m\in\mathbb{N}_{0};$ allowing $m=0$
reproduces the classical Laguerre polynomials. We also assume that
\begin{equation}
\alpha>m-1 \label{Type II parameter condition}%
\end{equation}
throughout this section. We refer the reader to the contributions
\cite{GUKM10, KMUG3, GUMM-Asymptotics} for full details of the Type II
exceptional $X_{m}$-Laguerre polynomials. We briefly discuss some of their
properties, and later, develop the spectral theory for the Type II exceptional
$X_{m}$-Laguerre differential expression.

\subsection{Properties of the Type II Exceptional $X_{m}$-Laguerre
Polynomials}

Choosing the factorization function $\phi_{2}(x),$ as given in (\ref{eq:seed5}%
), and letting $xL_{m}^{-\alpha}(x)$ be the factorization gauge, the classical
Laguerre differential expression (\ref{Classical Laguerre}) may be written as
\begin{align*}
-\ell^{\alpha}  &  =B_{m}^{\text{II,}\alpha}\circ A_{m}^{\text{II,}\alpha
}+\alpha-m,\quad\text{where\medskip}\\
A_{m}^{\text{II,}\alpha}[y](x)  &  =xL_{m}^{-\alpha}(x)y^{\prime}%
(x)+(\alpha-m)L_{m}^{-\alpha-1}(x)y(x)\quad\text{and\medskip}\\
B_{m}^{\text{II,}\alpha}[y](x)  &  =\frac{y^{\prime}(x)-y(x)}{L_{m}^{-\alpha
}(x)}\,.
\end{align*}

Based on this factorization, we define the Type II exceptional $X_{m}%
$-Laguerre expression $\ell_{m}^{\text{II,}\alpha}[\cdot]$ by
\begin{align}
\ell_{m}^{\text{II,}\alpha}[y]  &  =-\left(  A_{m}^{\text{II},\alpha+1}\circ
B_{m}^{\text{II},\alpha+1}[y]-m+\alpha+1\right)  \medskip\nonumber\\
&  =-\ell^{\alpha}[y](x)-2x(\log L_{m}^{-\alpha-1}(x))^{\prime}(y(x)-y^{\prime
}(x))+my(x)\medskip\nonumber\\
&  =-xy^{\prime\prime}(x)+\left(  -1-\alpha+x+2x\frac{(L_{m}^{-\alpha
-1}(x))^{\prime}}{L_{m}^{-\alpha-1}(x)}\right)  y^{\prime}(x)+\left(
m-2x\frac{(L_{m}^{-\alpha-1}(x))^{\prime}}{L_{m}^{-\alpha-1}(x)}\right)  y(x).
\label{Type II Expression}%
\end{align}

The Type II exceptional $X_{m}$-Laguerre polynomial $y=L_{m,n}^{\text{II}%
,\alpha}(x),$ where $n\geq m,$ satisfies the second-order differential
equation
\[
\ell_{m}^{\text{II,}\alpha}[y]=\lambda_{n}y\quad(0<x<\infty)
\]
where
\[
\lambda_{n}=n-m\quad(n\geq m)\,.
\]
Note that $\{\lambda_{n}\}_{n=m}^{\infty}=\mathbb{N}_{0}.$

The $n^{th}$ degree Type II exceptional $X_{m}$-Laguerre polynomial is
explicitly given by
\begin{align*}
L_{m,n}^{\text{II},\alpha}(x)  &  =-A_{m}^{\text{II,}\alpha+1}[L_{n-m}%
^{\alpha+1}](x)\\
&  =xL_{m}^{-\alpha-1}(x)L_{n-m-1}^{\alpha+2}(x)+(m-\alpha-1)L_{m}^{-\alpha
-2}(x)L_{n-m}^{\alpha+1}(x)\quad(n\geq m)\,.
\end{align*}

The sequence $\{L_{m,n}^{\text{II},\alpha}{}\}_{n=m}^{\infty}$ of Type II
exceptional $X_{m}$-Laguerre polynomials is orthogonal on $(0,\infty)$ with
respect to the weight function
\[
W_{m}^{\text{II},\alpha}(x)=\frac{x^{\alpha}e^{-x}}{(L_{m}^{-\alpha-1}%
(x))^{2}}\quad(x\in(0,\infty))\,.
\]

\begin{remark}
Requiring $\alpha>m-1$ is equivalent to $L_{m}^{-\alpha-1}(x)$ having no zeros
in $\left[  0,\infty\right)  $; see \cite[Proposition 4.1]{Odake-Sasaki}.
Notice that the function $1/(L_{m}^{-\alpha-1}(-x))^{2}$ is bounded and
bounded away from zero on $(0,\infty);$ hence all moments for the weight
function $W_{m}^{\textnormal{II},\alpha}$ on the interval $(0,\infty)$ exist and
are finite.
\end{remark}

In fact, in \cite{KMUG3}, the authors show that $\left\{  L_{m,n}%
^{\text{II},\alpha}{}\right\}  _{n=m}^{\infty}$ forms a complete orthogonal
set in the Hilbert space $L^{2}((0,\infty);W_{m}^{\text{II},\alpha}).$ With
$\left\Vert \cdot\right\Vert _{m}^{\text{II},\alpha}$ denoting the norm in
$L^{2}((0,\infty);W_{m}^{\text{II},\alpha})$ defined by%
\[
\left\Vert f\right\Vert _{m}^{\text{II},\alpha}=\left(  \int_{0}^{\infty
}\left\vert f(x)\right\vert ^{2}W_{m}^{\text{II},\alpha}(x)dx\right)
^{1/2}\quad(f\in L^{2}((0,\infty);W_{m}^{\text{II},\alpha}))
\]
in $L^{2}((0,\infty);W_{m}^{\text{II},\alpha}),$ derived from the inner
product%
\[
(f,g)_{m}^{\text{II},\alpha}:=\int_{0}^{\infty}f(x)\overline{g}(x)W_{m}%
^{\text{II},\alpha}(x)dx\quad(f,g\in L^{2}((0,\infty);W_{m}^{\text{II},\alpha
})),
\]
the norms of the Type II exceptional $X_{m}$-Laguerre polynomials are
explicitly given by%
\[
\left(  \left\Vert L_{m,n}^{\text{II},\alpha}\right\Vert _{m}^{\text{II}%
,\alpha}\right)  ^{2}=\frac{(\alpha+1+n-2m)\Gamma(\alpha+2+n-m)}{(n-m)!}%
\quad(n\geq m);
\]
see \cite{GUKM10} .

In \cite{GUMM-Asymptotics}, the authors establish the following two theorems
concerning properties of the zeros of $\left\{  L_{m,n}^{\text{II},\alpha}%
{}\right\}  _{n=m}^{\infty}.$

\begin{theorem}
\cite[Propositions 4.3, 4.4, and 4.5]{GUMM-Asymptotics} For $n\geq m,$ the
Type II exceptional $X_{m}$-Laguerre polynomial $L_{m,n}^{\textnormal{II},\alpha
}(x)$ has $n-m$ simple, positive zeros in $(0,\infty).$ Moreover,
$L_{m,n}^{\textnormal{II},\alpha}(x)$ has either $1$ or $0$ negative roots
according to, respectively, whether $m$ is odd or even.
\end{theorem}

\begin{theorem}
\cite[Corollary 4.1 and Proposition 4.8]{GUMM-Asymptotics} Let $\{j_{\alpha
,i}\}$ denote the sequence of positive zeros of the Bessel function of the
first kind $J_{\alpha}(x)$ listed in increasing order and let $\{x_{n,i}%
^{\alpha}\}_{i=1}^{n-m}$ denote the positive zeros of $L_{m,n}%
^{\textnormal{II},\alpha}(x)$ arranged in increasing order. Then%
\[
\lim_{n\rightarrow\infty}nx_{n,i}^{\alpha}=\frac{j_{\alpha,i}^{2}}{4}%
\quad(i\in\mathbb{N}).
\]
Furthermore, as $n\rightarrow\infty,$ the negative and complex roots of
$L_{m,n}^{\textnormal{II},\alpha}(z)$ converge to the zeros of $L_{m}^{-\alpha
-1}(z).$
\end{theorem}

\subsection{Type II Exceptional $X_{m}$-Laguerre Spectral Analysis}

In Lagrangian symmetric form, the Type II exceptional $X_{m}$- Laguerre
differential expression (\ref{Type II Expression}) is given by
\begin{align}
\ell_{m}^{\text{II,}\alpha}[y](x)  &  =\frac{1}{W_{m}^{\text{II},\alpha}%
(x)}\left(  -\left(  \frac{x^{\alpha+1}e^{-x}}{(L_{m}^{-\alpha-1}(x))^{2}%
}y^{\prime}(x)\right)  ^{\prime}\right. \label{Type II Symmetric Form}\\
&  +\left.  \left(  -\frac{mx^{\alpha}e^{-x}}{(L_{m}^{-\alpha-1}(x))^{2}%
}-\frac{2x^{\alpha+1}e^{-x}\left(  L_{m}^{-\alpha-1}(x)\right)  ^{\prime}%
}{\left(  L_{m}^{-\alpha-1}(x)\right)  ^{3}}\right)  y(x)\right)  .\nonumber
\end{align}

The maximal domain associated with $\ell_{m}^{\text{II,}\alpha}[\cdot]$ in the
Hilbert space $L^{2}\left(  (0,\infty);W_{m}^{II,\alpha}\right)  $ is defined
by
\[
\Delta_{m}^{\text{II,}\alpha}:=\left\{  f:(0,\infty)\rightarrow\mathbb{C}\mid
f,f^{\prime}\in AC_{\text{\textrm{loc}}}(0,\infty);f,\ell_{m}^{\text{II,}%
\alpha}[f]\in L^{2}((0,\infty);W_{m}^{\text{II},\alpha})\right\}  .
\]
The associated maximal operator
\[
T_{1,m}^{\text{II,}\alpha}:\mathcal{D}(T_{1,m}^{\text{II,}\alpha})\subset
L^{2}((0,\infty);W_{m}^{\text{II},\alpha})\rightarrow L^{2}((0,\infty
);W_{m}^{\text{II},\alpha}),
\]
is defined to be
\begin{align*}
T_{1,m}^{\text{II,}\alpha}f  &  =\ell_{m}^{\text{II,}\alpha}[f]\\
f\in\mathcal{D}(T_{1,m}^{\text{II,}\alpha}):  &  =\Delta_{m}^{\text{II,}%
\alpha}.
\end{align*}

For $f,g\in\Delta_{m}^{\text{II,}\alpha},$ Green's formula is
\[
\int_{0}^{\infty}\ell_{m}^{\text{II,}\alpha}[f](x)\overline{g}(x)W_{m}%
^{\text{II},\alpha}(x)dx=[f,g]_{m}^{\text{II,}\alpha}(x)\mid_{x=0}^{x=\infty
}+\int_{0}^{\infty}f(x)\ell_{m}^{\text{II,}\alpha}[\overline{g}](x)W_{m}%
^{\text{II},\alpha}(x)dx\,,
\]
where $[\cdot,\cdot]_{m}^{\text{II},\alpha}(\cdot)$ is the sesquilinear form
defined by
\begin{equation}
\left[  f,g\right]  _{m}^{\text{II},\alpha}(x):=\frac{x^{\alpha+1}e^{-x}%
}{(L_{m}^{-\alpha-1}(x))^{2}}(f(x)\bar{g}^{\prime}(x)-f^{\prime}(x)\bar
{g}(x))\quad(0<x<\infty), \label{Type II Sesquilinear Form}%
\end{equation}
and where
\[
\left[  f,g\right]  _{m}^{\text{II},\alpha}(x)\mid_{x=0}^{x=\infty}:=\left[
f,g\right]  _{m}^{\text{II},\alpha}(\infty)-\left[  f,g\right]  _{m}%
^{\text{II},\alpha}(0)\,.
\]
By Green's formula and the definition of $\Delta_{m}^{\text{II,}\alpha},$ both
limits
\[
\lbrack f,g]_{m}^{\text{II,}\alpha}(\infty):=\lim_{x\rightarrow\infty
}[f,g]_{m}^{\text{II,}\alpha}(x)\text{ and }[f,g]_{m}^{\text{II,}\alpha
}(0):=\lim_{x\rightarrow0^{+}}[f,g]_{m}^{\text{II,}\alpha}(x)
\]
exist and are finite for all $f,g\in\Delta_{m}^{\text{II,}\alpha}.$

The adjoint of the maximal operator in $L^{2}\left(  (0,\infty);W_{m}%
^{\text{II},\alpha}\right)  $ is the minimal operator $T_{0,m}^{\text{II}%
,\alpha},$ defined in $L^{2}\left(  (0,\infty);W_{m}^{\text{II},\alpha
}\right)  ,$ by%
\begin{align}
T_{0,m}^{\text{II},\alpha}f  &  =\ell_{m}^{\text{II,}\alpha}%
[f]\label{Type II Minimal Operator}\\
f\in\mathcal{D}(T_{0,m}^{\text{II},\alpha}):  &  =\{f\in\Delta_{m}%
^{\text{II,}\alpha}\mid\lbrack f,g]_{m}^{\text{II,}\alpha}\mid_{x=0}%
^{x=\infty}=0\text{ for all }g\in\Delta_{m}^{\text{II,}\alpha}\}.\nonumber
\end{align}

In the same manner as in the Type I exceptional $X_{m}$-Laguerre case, we seek
to find the self-adjoint extension $T_{m}^{\text{II},\alpha}$ in $L^{2}\left(
(0,\infty);W_{m}^{\text{II},\alpha}\right)  ,$ generated by $\ell
_{m}^{\text{II,}\alpha}[\cdot],$ which has the Type II exceptional $X_{m}%
$-Laguerre polynomials $\left\{  L_{m,n}^{\text{II},\alpha}\right\}
_{n=m}^{\infty}$ as eigenfunctions. As in the Type I case, we first need to
determine the deficiency index of the minimal operator $T_{0,m}^{\text{II}%
,\alpha}$ in $L^{2}\left(  (0,\infty);W_{m}^{\text{II},\alpha}\right)  .$ In
turn, this requires a study of the behavior of solutions near the singular
endpoints $x=0$ and $x=\infty$ of the differential expression
(\ref{Type II Expression}). This analysis is similar to the Type I case in the
previous section so we omit many of the details.

The point $x=0$ is a regular singular endpoint of (\ref{Type II Expression});
the Frobenius indicial equation is $r(r+\alpha)=0.$ Consequently, two linearly
independent solutions of $\ell_{m}^{\text{II,}\alpha}[y]=0$ on $(0,\infty)$
behave asymptotically like $z_{1}(x)=1$ and $z_{2}(x)=x^{-\alpha}$ near $x=0.$
Clearly $z_{1}\in L^{2}\left(  (0,\infty);W_{m}^{\text{II},\alpha}\right)  $
but a calculation shows that $z_{2}\in L^{2}\left(  (0,1);W_{m}^{\text{II}%
,\alpha}\right)  $ only when $\alpha<1.$ Consequently, $\ell_{m}%
^{\text{II,}\alpha}[\cdot]$ is in the limit-point case at $x=0$ when
$\alpha\geq1$ and is in the limit-circle case at $x=0$ when $\alpha<1.$ More
specifically, recalling (\ref{Type II parameter condition}),

\begin{enumerate}
\item[(i)] if $m=0$ (the classical Laguerre case), $\ell_{m}^{\text{II,}%
\alpha}[\cdot]$ is in the limit-circle case at $x=0$ when $-1<\alpha<1$ and in
the limit-point case when $\alpha\geq1;$

\item[(ii)] if $m=1$ (so, by (\ref{Type II parameter condition}), $\alpha>0),$
$\ell_{m}^{\text{II,}\alpha}[\cdot]$ is in the limit-circle case at $x=0$ when
$0<\alpha<1$ and in the limit-point case when $\alpha\geq1;$

\item[(iii)] if $m\geq2,$ then $\alpha>m-1\geq1$ and thus $\ell_{m}%
^{\text{II,}\alpha}[\cdot]$ is in the limit-point case at $x=0.$
\end{enumerate}

The point $x=\infty$ is an irregular singular endpoint of $\ell_{m}%
^{\text{II,}\alpha}[\cdot].$ Again, we can explicitly solve $\ell
_{m}^{\text{II,}\alpha}[y]=0$ on $(0,\infty)$ for a basis of solutions. One
solution is $y_{1}(x)=L_{m,m}^{\text{II},\alpha}(x):=L_{m}^{-\alpha-2}(x)$
which clearly belongs to $L^{2}\left(  (0,\infty);W_{m}^{\text{II},\alpha
}\right)  .$ A second solution $y_{2}(x)$ can be found by the reduction of
order method; this method shows that
\[
y_{2}(x)=L_{m}^{-\alpha-2}(x)\int_{a}^{x}\frac{e^{t}}{t^{\alpha+1}}\left(
\frac{L_{m}^{-\alpha-1}(t)}{L_{m}^{-\alpha-2}(t)}\right)  ^{2}dt\quad(x>0)
\]
where $a>0$ is fixed but otherwise arbitrary. An analysis similar to that
given in part (b) of Theorem \ref{LP/LC Type I Theorem} shows that
$y_{2}\notin L^{2}\left(  (x^{\ast},\infty);W_{m}^{\text{II},\alpha}\right)  $
for some $x^{\ast}>0.$ Consequently, $\ell_{m}^{\text{II,}\alpha}[\cdot]$ is
in the limit-point case at $x=\infty$ for any choice of $\alpha>m-1$.

When $\ell_{m}^{\text{II,}\alpha}[\cdot]$ is in the limit-circle case at
$x=0,$ the Glazman-Krein-Naimark theory requires that one appropriate boundary
condition be imposed in order to generate a self-adjoint extension of the
minimal operator $T_{0,m}^{\text{II},\alpha}$ in $L^{2}\left(  (0,\infty
);W_{m}^{\text{II},\alpha}\right)  .$ We are interested in a particular
self-adjoint extension, namely that operator $T_{m}^{\text{II,}\alpha}$ that
has the Type II exceptional $X_{m}$-Laguerre polynomials $\left\{
L_{m,n}^{\text{II,}\alpha}\right\}  _{n=m}^{\infty}$ as eigenfunctions. As in
the case of the Type I exceptional $X_{m}$-Laguerre case, we can take this
boundary condition to be%
\[
\lbrack f,1]_{m}^{\text{II,}\alpha}(0)=0,
\]
where $[\cdot,\cdot]_{m}^{\text{II,}\alpha}$ is the sesquilinear form given in
(\ref{Type II Sesquilinear Form}). This boundary condition simplifies to
\[
\lim_{x\rightarrow0^{+}}x^{\alpha+1}f^{\prime}(x)=0.
\]
We summarize this discussion in the following theorem.

\begin{theorem}
Let $m\in\mathbb{N}_{0}$ and $\alpha>m-1.$ Let $T_{0,m}^{\textnormal{II},\alpha}$
be the minimal operator, defined in \eqref{Type II Minimal Operator} in
$L^{2}((0,\infty);W_{m}^{\textnormal{II},\alpha})$ generated by the Type II
exceptional $X_{m}$-Laguerre differential expression $\ell_{m}%
^{\textnormal{II,}\alpha}[\cdot]$ given in \eqref{Type II Expression} or \eqref{Type II Symmetric Form}.

\begin{enumerate}
\item[$(a)$] The deficiency index of $T_{0,m}^{\textnormal{II},\alpha}$ is

\begin{enumerate}
\item[$(i)$] $(1,1)$ when $m=0$ and $-1<\alpha<1,$ or when $m=1$ and
$0<\alpha<1;$

\item[$(ii)$] $(0,0)$ when $m=0$ and $\alpha\geq1,$ or when $m=1$ and
$\alpha\geq1,$ or when $m\geq2.$
\end{enumerate}

\item[$(b)$] The operator
\[
T_{m}^{\textnormal{II},\alpha}:\mathcal{D}(T_{m}^{\textnormal{II},\alpha})\subset
L^{2}((0,\infty);W_{m}^{\textnormal{II},\alpha})\rightarrow L^{2}((0,\infty
);W_{m}^{\textnormal{II},\alpha}\,),
\]
defined by
\begin{align*}
T_{m}^{\textnormal{II},\alpha}f  &  =\ell_{m}^{\textnormal{II,}\alpha}[f]\\
f  &  \in\mathcal{D}(T_{m}^{\textnormal{II},\alpha}),
\end{align*}
is self-adjoint. The domain of $T_{m}^{\textnormal{II},\alpha}$ is given by

\begin{enumerate}
\item[$(i)$]
\[
\mathcal{D}(T_{m}^{\textnormal{II},\alpha}):=\{f\in\Delta_{m}^{\textnormal{II,}\alpha
}\mid\lim_{x\rightarrow0^{+}}x^{\alpha+1}f^{\prime}(x)=0\}
\]
when the deficiency index of $T_{0,m}^{\textnormal{II},\alpha}$ is $(1,1),$ or by

\item[$(ii)$]
\[
\mathcal{D}(T_{m}^{\textnormal{II},\alpha}):=\Delta_{m}^{\textnormal{II,}\alpha}%
\]
when the deficiency index of $T_{0,m}^{\textnormal{II},\alpha}$ is $(0,0).$
\end{enumerate}
\end{enumerate}

\noindent Moreover, in either case, $T_{m}^{\textnormal{II},\alpha}$ has the Type
II exceptional $X_{m}$-Laguerre polynomials $\{L_{m,n}^{\textnormal{II},\alpha
}\}_{n=m}^{\infty}$ as a complete set of eigenfunctions in $L^{2}%
((0,\infty);W_{m}^{\textnormal{II},\alpha}).$ Lastly, the spectrum of
$T_{m}^{\textnormal{II},\alpha}\ $consists only of eigenvalues and is given by
\[
\sigma(T_{m}^{\textnormal{II},\alpha})=\mathbb{N}_{0}\,.
\]

\end{theorem}

\section{A New Sequence of Exceptional Laguerre Polynomials: The Type III
Exceptional $X_{m}$-Laguerre Polynomials\label{Type III Section}}

The Type $\text{III}$ exceptional $X_{m}$-Laguerre polynomials
\[
\left\{  L_{m,n}^{\text{III},\alpha}\mid n=0,m+1,m+2,m+3,\ldots\right\}
\]
is a new class of exceptional $X_{m}$-Laguerre orthogonal polynomials for the
parameter range $-1<\alpha<0.$ They can be derived from the quasi-rational
eigenfunctions of the classical Laguerre differential expression
(\ref{Classical Laguerre}) and they can also be obtained from a transformation
of the Type I exceptional $X_{m}$-Laguerre differential expression
(\ref{Type I Expression}). Both of these derivations will be developed in
Section \ref{Type III from Type I}. In Section
\ref{Type III Laguerre Polynomials}, we introduce the Type III exceptional
$X_{m}$-Laguerre polynomials and derive several representations of them.
Section \ref{Norms} deals with the computation of the norms of these
polynomials in $L^{2}((0,\infty);W_{m}^{\text{III},\alpha})$. In Section
\ref{Completeness of Type III Laguerre Polynomials}, we prove that the
sequence of Type III exceptional\ $X_{m}$-Laguerre polynomials forms a
complete set of functions in $L^{2}((0,\infty);W_{m}^{\text{III},\alpha}).$ We
will deal with a comprehensive study of the location of the roots and an
asymptotic analysis of the roots for the Type III exceptional $X_{m}$-Laguerre
polynomials in Section \ref{Properties of Roots}. Lastly, in Section
\ref{Type III Spectral Analysis}, we construct a self-adjoint operator in
$L^{2}((0,\infty);W_{m}^{\text{III},\alpha}),$ generated by the second-order
Type III\ exceptional $X_{m}$-Laguerre differential expression, having the
sequence of Type III exceptional $X_{m}$-Laguerre polynomials as eigenfunctions.

\subsection{Two Derivations of the Type III Exceptional $X_{m}$-Laguerre
Differential Expression\label{Type III from Type I}}

Consider the transformation arising from the quasi-rational solution $\phi
_{2}(x)$ in (\ref{eq:seed5})
\begin{equation}
z(x)=x^{-\alpha}y(x)\,. \label{I to III transformation}%
\end{equation}
A calculation shows that
\begin{equation}
\ell_{m}^{\text{I},\alpha}[z](x)=x^{-\alpha}\ell_{m}^{\text{III},-\alpha
}[y](x)\,, \label{I to III transformation 2}%
\end{equation}
where
\begin{equation}
\ell_{m}^{\text{III,}\alpha}[y](x):=-xy^{\prime\prime}(x)+\left(
-1-\alpha+x+2x\frac{(L_{m}^{-\alpha-1}(-x))^{\prime}}{L_{m}^{-\alpha-1}%
(-x)}\right)  y^{\prime}(x)+(-m+\alpha)y(x) \label{Type III Expression}%
\end{equation}
or, equivalently, with the notation $M(g)(f(x)):=g(x)f(x),$%
\begin{equation}
M(x^{\alpha})\circ\ell_{m}^{\text{I,}\alpha}\circ M(x^{-\alpha})=\ell
_{m}^{\text{III},-\alpha}. \label{Gauge Transformation}%
\end{equation}
With regards to this identity, we say that the Type\ I and Type III
expressions are related by a \textit{gauge transformation. }

We call (\ref{Type III Expression}) the Type III exceptional $X_{m}$-Laguerre
differential expression. In Lagrangian symmetric form, this expression can be
written as
\begin{equation}
\ell_{m}^{\text{III,}\alpha}[y](x)=\frac{1}{W_{m}^{\text{III},\alpha}%
(x)}\left(  \left(  -\frac{x^{\alpha+1}e^{-x}}{(L_{m}^{-\alpha-1}(x))^{2}%
}y^{\prime}(x)\right)  ^{\prime}+\frac{(-m+\alpha)x^{\alpha}e^{-x}}%
{(L_{m}^{-\alpha-1}(x))^{2}}y(x)\right)  \,, \label{Type III Symmetric Form}%
\end{equation}
where
\[
W_{m}^{\text{III},\alpha}(x)=\frac{x^{\alpha}e^{-x}}{\left(  L_{m}^{-\alpha
-1}(-x)\right)  ^{2}}\quad(x\in(0,\infty)).
\]

\begin{remark}
As we will see below, the parameter range for the identity
\eqref{Gauge Transformation} is $0<\alpha<1.$ In this regard, we remark that
the Type III polynomials, which we show are solutions of
\[
\ell_{m}^{\textnormal{III,}\alpha}[y]=\lambda y,
\]
for a certain sequence of the eigenvalue parameter $\lambda,$ are related to
the L3 family of rational extensions of the isotonic oscillator which were
investigated by Grandati in \cite{grandati11}. From the point of view of
Schr\"{o}dinger operators, the parameter range $0<\alpha<1$ corresponds to a
potential with a weakly attracting singularity at the origin. Qualitatively,
this kind of singularity makes the physics of the system ambiguous and
requires the imposition of a boundary condition at the origin for a
well-defined eigenvalue problem; see Section \ref{Type III Spectral Analysis}.
\end{remark}

We note that all of the moments of $W_{m}^{\text{III},\alpha}$ exist and are
finite only when
\[
-1\,<\alpha<0.
\]

At this point, it is unclear if the eigenvalue problem
\begin{equation}
\ell_{m}^{\text{III,}\alpha}[y](x)=\lambda y(x) \label{Type III DE}%
\end{equation}
produces polynomial solutions for certain values of $\lambda.$ In the next
section, we will show that (\ref{Type III DE}) has polynomial solutions of
degrees $n=0$ and all $n\geq m+1.$ We now argue that there cannot be
polynomial solutions to (\ref{Type III DE}) of degrees $n=1,2,\ldots,m$ for
any value of $\lambda\in\mathbb{C}.$ Indeed, suppose $y=p(x)$ is a polynomial
solution to (\ref{Type III DE}) for some $\lambda\in\mathbb{C}.$ From
(\ref{Type III Expression}), it follows that the term%
\[
2x\frac{(L_{m}^{-\alpha-1}(-x))^{\prime}}{L_{m}^{-\alpha-1}(-x)}p^{\prime}(x)
\]
is a polynomial. However, since the roots of the Laguerre polynomial
$L_{m}^{-\alpha-1}(-x)$ are simple and negative, we see in fact that
$p^{\prime}(x)/L_{m}^{-\alpha-1}(-x)$ is a polynomial. Consequently, either
$p$ is a constant or a polynomial of degree $\geq m+1.$ More specifically, it
is the case, for some polynomial $q,$ that $p^{\prime}(x)=L_{m}^{-\alpha
-1}(-x)q(x);$ see Lemmas \ref{Lemma 1} and \ref{Lemma 2} below.

To see that (\ref{Type III DE}) has orthogonal polynomial eigenfunctions, we
turn to a special rational factorization of the classical Laguerre expression
(\ref{Classical Laguerre}). Indeed, the rational factorization function in
this case is $\phi_{3}(x),$ where $\phi_{3}$ is defined in (\ref{eq:seed6}),
and the corresponding gauge function is $xL_{m}^{-\alpha}(-x)$.

Define the first-order operators $A_{m}^{\text{III,}\alpha}$ and
$B_{m}^{\text{III,}\alpha}$ by
\begin{align*}
A_{m}^{\text{III,}\alpha}[y](x)  &  :=xL_{m}^{-\alpha}(-x)y^{\prime
}(x)-(m+1)L_{m+1}^{-\alpha-1}(-x)y(x)\\
B_{m}^{\text{III,}\alpha}[y](x)  &  :=\frac{y^{\prime}(x)}{L_{m}^{-\alpha
}(-x)}\,.
\end{align*}

\begin{lemma}
\label{Properties of First Order Operators}The operators $A_{m}%
^{\textnormal{III,}\alpha}$ and $B_{m}^{\textnormal{III,}\alpha}$ satisfy the
following factorization properties:

\begin{enumerate}
\item[$(a)$] $-\ell^{\alpha}=B_{m}^{\textnormal{III,}\alpha}\circ A_{m}%
^{\textnormal{III,}\alpha}+m+1,$ where $\ell^{\alpha}$ is the classical Laguerre
second-order differential expression defined in
\eqref{Classical Laguerre};\medskip

\item[$(b)$] $-\ell_{m}^{\textnormal{III,}\alpha}=A_{m}^{\textnormal{III,}\alpha
+1}\circ B_{m}^{\textnormal{III,}\alpha+1}+m-\alpha.$
\end{enumerate}

\begin{proof}
The proofs of these identities are similar so we give only the proof of part
$($b$)$. Our proof will make repeated use of two facts:%
\begin{equation}
(L_{n}^{\alpha}(-x))^{\prime}=L_{n-1}^{\alpha+1}(-x)\text{ for any }%
\alpha>-1\text{ and }n\in\mathbb{N}_{0}, \label{Key Identity 1}%
\end{equation}
and%
\begin{equation}
x\left(  L_{m+1}^{-\alpha-2}(-x)\right)  ^{\prime\prime}+(x-\alpha-1)\left(
L_{m+1}^{-\alpha-2}(-x)\right)  ^{\prime}-(m+1)L_{m+1}^{-\alpha-2}(-x)=0;
\label{Key Identity 2}%
\end{equation}
see \cite[Chapter V, (5.1.2) and (5.1.14)]{Szego}. Now
\begin{align}
A_{m}^{\text{III},\alpha+1}\left(  B_{m}^{\text{III,}\alpha+1}[y]\right)   &
=xL_{m}^{-\alpha-1}(-x)\left(  \frac{y^{\prime}}{L_{m}^{-\alpha-1}%
(-x)}\right)  ^{\prime}-(m+1)L_{m+1}^{-\alpha-2}(-x)\left(  \frac{y^{\prime}%
}{L_{m}^{-\alpha-1}(-x)}\right) \nonumber\\
&  =xL_{m}^{-\alpha-1}(-x)\left(  \frac{L_{m}^{-\alpha-1}(-x)y^{\prime\prime
}-L_{m-1}^{-\alpha}(-x)y^{\prime}}{\left(  L_{m}^{-\alpha-1}(-x)\right)  ^{2}%
}\right)  -(m+1)\frac{L_{m+1}^{-\alpha-2}(-x)}{L_{m}^{-\alpha-1}(-x)}%
y^{\prime}\nonumber\\
&  =xy^{\prime\prime}+\left(  \frac{-xL_{m-1}^{-\alpha}(-x)-(m+1)L_{m+1}%
^{-\alpha-2}(-x)}{L_{m}^{-\alpha-1}(-x)}\right)  y^{\prime}.
\label{Partial Identity}%
\end{align}
Moreover,
\begin{align}
&  \frac{-xL_{m-1}^{-\alpha}(-x)-(m+1)L_{m+1}^{-\alpha-2}(-x)}{L_{m}%
^{-\alpha-1}(-x)}\nonumber\\
&  =\frac{-x\left(  L_{m+1}^{-\alpha-2}(-x)\right)  ^{\prime\prime
}-(m+1)L_{m+1}^{-\alpha-2}(-x)}{L_{m}^{-\alpha-1}(-x)}\text{ by }%
(\text{\ref{Key Identity 1}})\nonumber\\
&  =\frac{(x-\alpha-1)\left(  L_{m+1}^{-\alpha-2}(-x)\right)  ^{\prime
}-2(m+1)L_{m+1}^{-\alpha-2}(-x)}{L_{m}^{-\alpha-1}(-x)}\text{ from
}(\text{\ref{Key Identity 2}})\nonumber\\
&  =\frac{(x-\alpha-1)L_{m}^{-\alpha-1}(-x)-2(m+1)L_{m+1}^{-\alpha-2}%
(-x)}{L_{m}^{-\alpha-1}(-x)}\text{ by }(\text{\ref{Key Identity 1}%
})\nonumber\\
&  =\frac{(x-\alpha-1)L_{m}^{-\alpha-1}(-x)-2x\left(  L_{m+1}^{-\alpha
-2}(-x)\right)  ^{\prime\prime}+(2\alpha+2-2x)\left(  L_{m+1}^{-\alpha
-2}(-x)\right)  ^{\prime}}{L_{m}^{-\alpha-1}(-x)}\text{ by }%
(\text{\ref{Key Identity 2}})\nonumber\\
&  =\frac{(x-\alpha-1)L_{m}^{-\alpha-1}(-x)-2x\left(  L_{m}^{-\alpha
-1}(-x)\right)  ^{\prime}+(2\alpha+2-2x)L_{m}^{-\alpha-1}(-x)}{L_{m}%
^{-\alpha-1}(-x)}\text{ from }(\text{\ref{Key Identity 1}}).\nonumber\\
&  =1+\alpha-x-2x\frac{\left(  L_{m}^{-\alpha-1}(-x)\right)  ^{\prime}}%
{L_{m}^{-\alpha-1}(-x)}. \label{Partial Identity 2}%
\end{align}
Substitution of (\ref{Partial Identity 2}) into (\ref{Partial Identity})
yields%
\[
A_{m}^{\text{III},\alpha+1}\left(  B_{m}^{\text{III,}\alpha+1}[y]\right)
=xy^{\prime\prime}+\left(  1+\alpha-x-2x\frac{\left(  L_{m}^{-\alpha
-1}(-x)\right)  ^{\prime}}{L_{m}^{-\alpha-1}(-x)}\right)  y^{\prime};
\]
adding the term $(m-\alpha)y$ to both sides of this latter identity completes
the proof.
\end{proof}
\end{lemma}

\begin{remark}
With reference to \eqref{Rational Factorization} and
\eqref{Darboux Transformation}, where the reader will notice that the
parameters $\lambda$ in both expressions are equal, we could define the Type
III exceptional $X_{m}$-Laguerre differential expression by
\begin{align*}
\widetilde{\ell}_{m}^{\textnormal{III,}\alpha}[y](x):  &  =-\left(  A_{m}%
^{\textnormal{III,}\alpha+1}\circ B_{m}^{\textnormal{III,}\alpha+1}+m+1\right)
[y](x)\\
&  =-xy^{\prime\prime}(x)+\left(  -1-\alpha+x+2x\frac{(L_{m}^{-\alpha
-1}(-x))^{\prime}}{L_{m}^{-\alpha-1}(-x)}\right)  y^{\prime}(x)+(-m-1)y(x).
\end{align*}
In this case, we would not have the identity \eqref{Gauge Transformation};
however, by mimicking the proof of Theorem
\ref{Polynomial Solutions to Type III Expression} below, the Type III
exceptional $X_{m}$-Laguerre polynomial $y=L_{m,n}^{\textnormal{III},\alpha}(x)$
can be shown to be a solution of the eigenvalue equation
\[
\widetilde{\ell}_{m}^{\textnormal{III,}\alpha}[y](x)=(n-m-1)y(x)\quad
(n=0,m+1,m+2,m+3,\ldots).
\]

\end{remark}

\subsection{The Type $\text{III}$ Exceptional $X_{m}$-Laguerre
Polynomials\label{Type III Laguerre Polynomials}}

For the remainder of this section, we assume that $-1<\alpha<0.$

Similar to how we introduced the Type I and Type II exceptional $X_{m}%
$-Laguerre polynomials, we define the $n^{th}$ degree Type III exceptional
$X_{m}$-Laguerre polynomial by%
\begin{equation}
L_{m,n}^{\text{III},\alpha}(x):=\left\{
\begin{array}
[c]{ll}%
-A_{m}^{\text{III},\alpha+1}[L_{n-m-1}^{\alpha+1}](x) & \text{if }n\geq m+1\\
1 & \text{if }n=0.
\end{array}
\right.  \label{Definition of Type III Laguerre}%
\end{equation}
From the definition of $A_{m}^{\text{III},\alpha+1}$, a calculation shows
that
\begin{equation}
L_{m,n}^{\text{III},\alpha}(x)=%
\begin{cases}
xL_{n-m-2}^{\alpha+2}(x)L_{m}^{-\alpha-1}(-x)+(m+1)L_{n-m-1}^{\alpha
+1}(x)L_{m+1}^{-\alpha-2}(-x), & \text{if }n\geq m+1\\
1 & \text{if }n=0.
\end{cases}
\label{Type III Laguerre rep}%
\end{equation}

The following lemmas (Lemmas \ref{Lemma 1} and \ref{Lemma 2}) are critical for
several reasons. Indeed, they will ultimately help show that $y=L_{m,n}%
^{III,\alpha}(x)$ is a solution of the eigenvalue equation
\begin{equation}
\ell_{m}^{\text{III,}\alpha}[y](x)=\lambda_{n}y(x),
\label{L3 eigenvalue equation}%
\end{equation}
where
\begin{equation}
\lambda_{n}=n-m+\alpha\quad(n=0,m+1,m+2,m+3,\ldots). \label{L3 eigenvalues}%
\end{equation}
In addition, both lemmas give new characterizations of the Type III
exceptional $X_{m}$-Laguerre polynomials and lead to an additional
representation (Theorem
\ref{New Representation of Type III Laguerre Polynomials}) of these
polynomials. Lastly, these lemmas will be critically important in our analysis
of the location of the roots (Lemma \ref{Lemma Representation} and Theorem
\ref{Theorem Interlacing}) of $\left\{  L_{m,n}^{\text{III},\alpha}\right\}  $
and to proving root asymptotic results (Theorem \ref{Theorem Asymptotic}) of
these roots.

\begin{lemma}
\label{Lemma 1} For $k\in\mathbb{N},$
\begin{equation}
\frac{\left(  L_{m,m+k}^{\textnormal{III},\alpha}(x)\right)  ^{\prime}}%
{L_{m}^{-\alpha-1}(-x)}=-xL_{k-3}^{\alpha+3}(x)+(\alpha+2-x)L_{k-2}^{\alpha
+2}(x)+(m+1)L_{k-1}^{\alpha+1}(x). \label{Key Relationship}%
\end{equation}

\begin{proof}
Recall the representation $($\ref{Type III Laguerre rep}$)$:
\[
L_{m,m+k}^{\text{III},\alpha}(x)=xL_{k-2}^{\alpha+2}(x)L_{m}^{-\alpha
-1}(-x)+(m+1)L_{k-1}^{\alpha+1}(x)L_{m+1}^{-\alpha-2}(-x).
\]
Using the Laguerre identity
\begin{equation}
(L_{k}^{\alpha}(x))^{\prime}=-L_{k-1}^{\alpha+1}(x), \label{Lemma 1 - 1}%
\end{equation}
we see that%
\begin{align*}
\left(  L_{m,m+k}^{\text{III},\alpha}(x)\right)  ^{\prime}  &  =L_{k-2}%
^{\alpha+2}(x)L_{m}^{-\alpha-1}(-x)-xL_{k-3}^{\alpha+3}(x)L_{m}^{-\alpha
-1}(-x)\\
&  +xL_{k-2}^{\alpha+2}(x)L_{m-1}^{-\alpha}(-x)-(m+1)L_{k-2}^{\alpha
+2}(x)L_{m+1}^{-\alpha-2}(-x)\\
&  +(m+1)L_{k-1}^{\alpha+1}(x)L_{m}^{-\alpha-1}(-x).
\end{align*}
Therefore,%
\begin{align}
\frac{\left(  L_{m,m+k}^{\text{III},\alpha}(x)\right)  ^{\prime}}%
{L_{m}^{-\alpha-1}(-x)}=  &  -xL_{k-3}^{\alpha+3}(x)+(m+1)L_{k-1}^{\alpha
+1}(x)\nonumber\\
\quad\quad &  +L_{k-2}^{\alpha+2}(x)\left(  1+x\frac{L_{m-1}^{-\alpha}%
(-x)}{L_{m}^{-\alpha-1}(-x)}-(m+1)\frac{L_{m+1}^{-\alpha-2}(-x)}%
{L_{m}^{-\alpha-1}(-x)}\right)  . \label{Lemma 1 - 2}%
\end{align}
Now, using \eqref{Lemma 1 - 1},
\begin{align}
&  1+x\frac{L_{m-1}^{-\alpha}(-x)}{L_{m}^{-\alpha-1}(-x)}-(m+1)\frac
{L_{m+1}^{-\alpha-2}(-x)}{L_{m}^{-\alpha-1}(-x)}\label{Lemma 1 -3}\\
&  =\frac{L_{m}^{-\alpha-1}(-x)+xL_{m-1}^{-\alpha}(-x)-(m+1)L_{m+1}%
^{-\alpha-2}(-x)}{L_{m}^{-\alpha-1}(-x)}\nonumber\\
&  =\frac{x(L_{m+1}^{-\alpha-2}(-x))^{\prime\prime}+(L_{m+1}^{-\alpha
-2}(-x))^{\prime}-(m+1)L_{m+1}^{-\alpha-2}(-x)}{L_{m}^{-\alpha-1}%
(-x)}.\nonumber
\end{align}
Since $y=L_{m+1}^{-\alpha-2}(x)$ satisfies $xy^{\prime\prime}+(-\alpha
-1-x)y^{\prime}+(m+1)y=0,$ a simple calculation shows that $y=L_{m+1}%
^{-\alpha-2}(-x)$ satisfies%
\[
xy^{\prime\prime}+(-\alpha-1+x)y^{\prime}-(m+1)y=0.
\]
Hence \eqref{Lemma 1 -3} becomes%
\begin{align}
&  1+x\frac{L_{m-1}^{-\alpha}(-x)}{L_{m}^{-\alpha-1}(-x)}-(m+1)\frac
{L_{m+1}^{-\alpha-2}(-x)}{L_{m}^{-\alpha-1}(-x)}\nonumber\\
&  =\frac{(\alpha+2-x)(L_{m+1}^{-\alpha-2}(-x))^{\prime}}{L_{m}^{-\alpha
-1}(-x)}=\alpha+2-x \label{Lemma 1 -4}%
\end{align}
since $(L_{m+1}^{-\alpha-2}(-x))^{\prime}=L_{m}^{-\alpha-1}(-x).$ Substituting
\eqref{Lemma 1 -4} into \eqref{Lemma 1 - 2} establishes \eqref{Key Relationship}.
\end{proof}
\end{lemma}

\begin{lemma}
\label{Lemma 2} For $k\in\mathbb{N},$
\begin{equation}
\left(  L_{m,m+k}^{\textnormal{III},\alpha}(x)\right)  ^{\prime}=(m+k)L_{k-1}%
^{\alpha+1}(x)L_{m}^{-\alpha-1}(-x). \label{Lemma 2 Claim}%
\end{equation}
In particular, for all $k\in\mathbb{N}$, the Type III exceptional $X_{m}%
$-Laguerre polynomial $L_{m,m+k}^{\textnormal{III},\alpha}(x)$ has a local
extremum at each of the $m$ roots of $L_{m}^{-\alpha-1}(-x)$ and at each of
the $k-1$ roots of $L_{k-1}^{\alpha+1}(x)$.

\begin{proof}
The Laguerre polynomial $y=L_{k-1}^{\alpha+1}(x)$ is a solution of Laguerre's
equation%
\[
xy^{\prime\prime}+(\alpha+2-x)y^{\prime}+(k-1)y=0.
\]
Consequently, we see that the right-hand side of \eqref{Key Relationship}
simplifies to%
\begin{align*}
&  -xL_{k-3}^{\alpha+3}(x)+(\alpha+2-x)L_{k-2}^{\alpha+2}(x)+(m+1)L_{k-1}%
^{\alpha+1}(x)\\
&  =-x(L_{k-1}^{\alpha+1}(x))^{\prime\prime}-(\alpha+2-x)(L_{k-1}^{\alpha
+1}(x))^{\prime}+(m+1)L_{k-1}^{\alpha+1}(x)\\
&  =(m+k)L_{k-1}^{\alpha+1}(x).
\end{align*}
The result now follows from this identity and \eqref{Key Relationship}.

A degree count implies that all roots of $\left(  L_{m,m+k}^{\text{III}
,\alpha}(x)\right)  ^{\prime}$ must be simple.
\end{proof}
\end{lemma}

We are now in position to prove the following theorem which will establish
(\ref{L3 eigenvalue equation}) and (\ref{L3 eigenvalues}).

\begin{theorem}
\label{Polynomial Solutions to Type III Expression}For $n=0,$ $m+1,$ $m+2,$
$m+3,\ldots,$ the function $y=L_{m,n}^{\textnormal{III},\alpha}(x)$ is a solution
of \eqref{L3 eigenvalue equation}, where $\lambda_{n}$ is given in \eqref{L3 eigenvalues}.

\begin{proof}
The proof is straightforward when $n=0$ so we assume $n\geq m+1.$ With
$y=L_{m,n}^{\text{III},\alpha}(x),$ we see from Lemma
\ref{Properties of First Order Operators} $($b$)$ that%
\begin{align*}
\ell_{m}^{\text{III},\alpha}[y](x)  &  =-A_{m}^{\text{III},\alpha+1}%
B_{m}^{\text{III},\alpha+1}[L_{m,n}^{\text{III},\alpha}](x)+(-m+\alpha
)L_{m,n}^{\text{III},\alpha}(x)\\
&  =-A_{m}^{\text{III},\alpha+1}\left[  (L_{m,n}^{\text{III},\alpha
}(x))^{\prime}/L_{m}^{-\alpha-1}(-x)\right]  +(-m+\alpha)L_{m,n}%
^{\text{III},\alpha}(x)\text{ by definition of }B_{m}^{\text{III},\alpha+1}\\
&  =-A_{m}^{\text{III},\alpha+1}[nL_{n-m-1}^{\alpha+1}](x)+(-m+\alpha
)L_{m,n}^{\text{III},\alpha}(x)\text{ by Lemma \ref{Lemma 2}}\\
&  =nL_{m,n}^{\text{III},\alpha+1}(x)+(-m+\alpha)L_{m,n}^{\text{III},\alpha
}(x)\text{ by }(\text{\ref{Definition of Type III Laguerre}})\\
&  =(n-m+\alpha)y(x).
\end{align*}

\end{proof}
\end{theorem}

The next two results give new representations of the Type III exceptional
$X_{m}$-Laguerre polynomials.

\begin{theorem}
\label{New Representation of Type III Laguerre Polynomials}For $k\in
\mathbb{N},$%
\begin{equation}
L_{m,m+k}^{\textnormal{III},\alpha}(x)=(m+k)\int_{0}^{x}L_{k-1}^{\alpha+1}%
(t)L_{m}^{-\alpha-1}(-t)dt+(m+1)\binom{k+\alpha}{k-1}\dbinom{m-\alpha-1}{m+1}.
\label{Simplest Representation}%
\end{equation}

\begin{proof}
This follows immediately from \eqref{Lemma 2 Claim} and the normalizations
\[
L_{k}^{\alpha}(0)=\binom{k+\alpha}{k}\quad\text{and}\quad L_{m,m+k}%
^{\text{III},\alpha}(0)=(m+1)L_{k-1}^{\alpha+1}(0)L_{m+1}^{-\alpha-2}(0).
\]

\end{proof}
\end{theorem}

The following representation of the Type $\text{III}$ exceptional $X_{m}%
$-Laguerre polynomials will be important for determining the location of their zeros.

\begin{lemma}
\label{Lemma Representation} For $m,k\in\mathbb{N}$ and $-1<\alpha<0,$ the
Type $\text{III}$ exceptional $X_{m}$-Laguerre polynomial $L_{m,m+k}%
^{\textnormal{III},\alpha}(x)$ can be written as
\begin{align}
L_{m,m+k}^{\textnormal{III},\alpha}(x)  &  =(k+\alpha)L_{k-2}^{\alpha+1}%
(x)L_{m}^{-\alpha-1}(-x)+(m+1)L_{k-1}^{\alpha+1}(x)L_{m+1}^{-\alpha
-1}(-x)\label{Eq Alternative}\\
&  \quad-(m+k)L_{k-1}^{\alpha+1}(x)L_{m}^{-\alpha-1}(-x).\nonumber
\end{align}

\begin{proof}
Recall \eqref{Type III Laguerre rep}:
\[
L_{m,m+k}^{\text{III},\alpha}(x)=xL_{k-2}^{\alpha+2}(x)L_{m}^{-\alpha
-1}(-x)+(m+1)L_{k-1}^{\alpha+1}(x)L_{m+1}^{-\alpha-2}(-x)\,.
\]
From \cite[p. 102, (5.1.14)]{Szego},
\[
x(L_{n}^{\alpha}(x))^{\prime}=-(n+\alpha)L_{n-1}^{\alpha}(x)+nL_{n}^{\alpha
}(x)\,,
\]
we see that
\begin{align}
xL_{k-2}^{\alpha+2}(x)  &  =-x\left(  L_{k-1}^{\alpha+1}(x)\right)  ^{\prime
}\nonumber\\
&  =(k+\alpha)L_{k-2}^{\alpha+1}(x)-(k-1)L_{k-1}^{\alpha+1}(x). \label{B}%
\end{align}

Likewise, from the identity $($see \cite[p. 102, (5.1.13)]{Szego}$)$
\[
L_{n}^{\alpha}(x)=L_{n}^{\alpha+1}(x)-L_{n-1}^{\alpha+1}(x)\,,
\]
we obtain
\begin{equation}
L_{m+1}^{-\alpha-2}(-x)=L_{m+1}^{-\alpha-1}(-x)-L_{m}^{-\alpha-1}(-x)\,.
\label{C}%
\end{equation}
Substituting \eqref{B} and \eqref{C} into \eqref{Type III Laguerre rep} yields
$($\ref{Eq Alternative}$)$.
\end{proof}
\end{lemma}

\begin{remark}
Our discussion to this point shows that if we take \eqref{Eq Alternative}
$($or \eqref{Simplest Representation}$)$ as our \emph{definition} of the Type
III exceptional $X_{m}$-Laguerre polynomials, they are orthogonal polynomials
for $-1<\alpha<0.$ Regardless of this parameter restriction, the polynomial
defined in either \eqref{Eq Alternative} $($or
\eqref{Simplest Representation}$)$ is of degree $m+k.$ A natural question to
ask is whether or not these polynomials are orthogonal, in some sense, for
values of $\alpha\notin(-1,0).$ It would be interesting to look into this
question even when some of the associated moments do not exist.
\end{remark}

We note that the Type $\text{III}$ exceptional $X_{m}$-Laguerre polynomials
are negative at the origin; that is to say, for $-1<\alpha<0$ and
$k\in\mathbb{N},$ we have
\begin{equation}
L_{m,m+k}^{\text{III},\alpha}(0)<0. \label{Negativity of Type III at 0}%
\end{equation}
To see this, recall from Theorem
\ref{New Representation of Type III Laguerre Polynomials} that
\[
L_{m,m+k}^{\text{III},\alpha}(0)=(m+1)L_{k-1}^{\alpha+1}(0)L_{m+1}^{-\alpha
-2}(0).
\]
Now, in general,%
\[
L_{n}^{\alpha}(0)=\dbinom{n+\alpha}{n}=\frac{\Gamma(n+\alpha+1)}{\Gamma
(\alpha+1)n!}=\frac{(1+\alpha)(2+\alpha)\ldots(n+\alpha)}{n!},
\]
so $L_{k-1}^{\alpha+1}(0)>0$ when $-1<\alpha<0.$ However, this parameter range
implies that $-1<-\alpha-1<0$ while $-\alpha+j>0$ for $j=0,1,\ldots m-1$;
thus
\begin{equation}
L_{m+1}^{-\alpha-2}(0)=\frac{(-\alpha-1)(-\alpha)(-\alpha+1)\cdots
(-\alpha+m-1)}{(m+1)!}<0. \label{Negative}%
\end{equation}
From (\ref{Negative}), the inequality in (\ref{Negativity of Type III at 0})
now follows. The negativity of $L_{m,m+k}^{\text{III},\alpha}(0)$ turns out to
be essential in our analysis (see Section \ref{Properties of Roots}) of
determining the location of the roots of $L_{m,m+k}^{\text{III},\alpha}(x)$.

\subsection{The Norms of the Type III Exceptional $X_{m}$-Laguerre
Polynomials\label{Norms}}

We now compute the norms of these Type III polynomials.

\begin{theorem}
\label{Norms of Type III Polynomials}Suppose $-1\,<\alpha<0.$ The Type\ III
exceptional $X_{m}$-Laguerre polynomials
\[
\left\{  L_{m,n}^{\textnormal{III},\alpha}\mid n=0,m+1,m+2,m+3,\ldots\right\}
\]
are orthogonal in the Hilbert space $L^{2}((0,\infty);W_{m}%
^{\textnormal{III},\alpha})$ and the norms are explicitly given by
\begin{equation}
\left(  \Vert L_{m,n}^{\textnormal{III},\alpha}\Vert_{m}^{\textnormal{III},\alpha
}\right)  ^{2}=\int_{0}^{\infty}(L_{m,n}^{\textnormal{III},\alpha}(x))^{2}%
W_{m}^{\textnormal{III},\alpha}(x)dx=%
\begin{cases}
\displaystyle\frac{n\,\Gamma(n-m+\alpha+1)}{(n-m-1)!} & \text{if }n\geq
m+1\\[12pt]%
\displaystyle\frac{\Gamma(\alpha+1)\Gamma(-\alpha)m!}{\Gamma(m-\alpha)} &
\text{if }n=0.
\end{cases}
\label{eq:Type3norm}%
\end{equation}

\begin{proof}
We compute the norms in this proof; the orthogonality will follow directly
from the self-adjointness of the operator $T_{m}^{\text{III},\alpha}$ in
Theorem \ref{Type III Self-Adjoint Operator} in Section
\ref{Type III Spectral Analysis} below. The proof, when $n\geq m+1,$ rests on
the following adjoint relationship for the $A_{m}^{\text{III,}\alpha}{}$ and
$B_{m}^{\text{III,}\alpha}{}$ operators
\begin{equation}
B_{m}^{\text{III,}\alpha}[f](x)\,g(x)W^{\alpha}(x)+A_{m}^{\text{III,}\alpha
}[g](x)f(x)W_{m}^{\text{III},\alpha-1}(x)=\frac{d}{dx}\left(  \frac{W^{\alpha
}(x)}{L_{m}^{-\alpha}(-x)}f(x)g(x)\right)  , \label{eq:A3B3adjoint}%
\end{equation}
where we take $f=f(x),g=g(x)$ to be polynomials and where $W^{\alpha
}(x)=x^{\alpha}e^{-x}$ is the classical Laguerre weight. To prove this, divide
the left-hand side of (\ref{eq:A3B3adjoint}) by
\[
\frac{W^{\alpha}(x)}{L_{m}^{-\alpha}(-x)}f(x)g(x)
\]
to obtain
\[
\frac{L_{m}^{-\alpha}(-x)B_{m}^{\text{III,}\alpha}[f](x)}{f(x)}+\frac
{A_{m}^{\text{III,}\alpha}[g](x)}{xL_{m}^{-\alpha}(-x)g(x)};
\]
on the other hand, a tedious calculation shows that
\begin{align*}
&  \frac{d}{dx}\log\left(  \frac{W^{\alpha}(x)}{L_{m}^{-\alpha}(-x)}%
f(x)g(x)\right)  \medskip\\
&  =\frac{f^{\prime}(x)}{f(x)}+\frac{xL_{m}^{-\alpha}(-x)g^{\prime}(x)+\left(
\alpha L_{m}^{-\alpha}(-x)-xL_{m}^{-\alpha}(-x)-x(L_{m}^{-\alpha}%
(-x))^{\prime}\right)  g(x)}{xL_{m}^{-\alpha}(-x)g(x)}\medskip\\
&  =\frac{f^{\prime}(x)}{f(x)}+\frac{xL_{m}^{-\alpha}(-x)g^{\prime}(x)+\left(
-x(L_{m+1}^{-\alpha-1}(-x))^{\prime\prime}+(\alpha-x)(L_{m+1}^{-\alpha
-1}(-x))^{\prime}\right)  g(x)}{xL_{m}^{-\alpha}(-x)g(x)}\medskip\\
&  =\frac{f^{\prime}(x)}{f(x)}+\frac{xL_{m}^{-\alpha}(-x)g^{\prime
}(x)-(m+1)L_{m+1}^{-\alpha-1}(-x)g(x)}{xL_{m}^{-\alpha}(-x)g(x)}\text{ from
}(\text{\ref{Key Identity 2}})\text{ with }\alpha\rightarrow\alpha-1\medskip\\
&  =\frac{L_{m}^{-\alpha}(-x)B_{m}^{\text{III,}\alpha}[f](x)}{f(x)}%
+\frac{A_{m}^{\text{III,}\alpha}[g](x)}{xL_{m}^{-\alpha}(-x)g(x)}\text{ by
definition of }A_{m}^{\text{III,}\alpha}\text{ and }B_{m}^{\text{III,}\alpha}.
\end{align*}
This establishes $($\ref{eq:A3B3adjoint}$).$ Setting $f=A_{m}^{\text{III,}%
\alpha}[L_{n}^{\alpha}]$ and $g=L_{n}^{\alpha}$ in $($\ref{eq:A3B3adjoint}$)$,
integrating and using Lemma \ref{Properties of First Order Operators}, Part
(a), and (\ref{Definition of Type III Laguerre}) gives
\begin{align}
-  &  (n+m+1)\int_{0}^{\infty}\left(  L_{n}^{\alpha}(x)\right)  ^{2}W^{\alpha
}(x)dx+\int_{0}^{\infty}\left(  L_{m,m+n+1}^{\text{III},\alpha-1}(x)\right)
^{2}W_{m}^{\text{III},\alpha-1}(x)dx\label{Norm relationship}\\
&  =-\left(  \frac{L_{m,m+n+1}^{\text{III,}\alpha-1}(x)L_{n}^{\alpha}%
(x)}{L_{m}^{-\alpha}(-x)}W^{\alpha}(x)\right)
\bigg|%
_{x=0}^{x=\infty}.\nonumber
\end{align}
In $($\ref{Norm relationship}$),$ we are assuming that $-1<\alpha-1<0$ so
$0<\alpha<1.$ It follows that
\[
\left(  \frac{L_{m,m+n+1}^{\text{III,}\alpha-1}(x)L_{n}^{\alpha}(x)}%
{L_{m}^{-\alpha}(-x)}W^{\alpha}(x)\right)
\bigg|%
_{x=0}^{x=\infty}=\left(  \frac{L_{m,m+n+1}^{\text{III,}\alpha-1}%
(x)L_{n}^{\alpha}(x)}{L_{m}^{-\alpha}(-x)}x^{\alpha}e^{-x}\right)
\bigg|%
_{x=0}^{x=\infty}=0.
\]
Hence, we see that
\[
\int_{0}^{\infty}\left(  L_{m,m+n+1}^{\text{III},\alpha-1}(x)\right)
^{2}W_{m}^{\text{III},\alpha-1}(x)dx=(n+m+1)\int_{0}^{\infty}\left(
L_{n}^{\alpha}(x)\right)  ^{2}W^{\alpha}(x)dx.
\]
Since $($see \cite[Chapter V, (5.1.1)]{Szego}$)$%
\[
\int_{0}^{\infty}\left(  L_{n}^{\alpha}(x)\right)  ^{2}W^{\alpha}%
(x)dx=\frac{\Gamma(n+\alpha+1)}{n!}\quad(j\in\mathbb{N}_{0}),
\]
we see from (\ref{Norm relationship}) that
\[
\int_{0}^{\infty}\left(  L_{m,m+n+1}^{\text{III},\alpha-1}(x)\right)
^{2}W_{m}^{\text{III},\alpha-1}(x)dx=\frac{(n+m+1)\Gamma(n+\alpha+1)}{n!}.
\]
Replacing $\alpha$ by $\alpha+1$ and $n+m+1$ by $n$ yields%
\[
\int_{0}^{\infty}(L_{m,n}^{\text{III},\alpha}(x))^{2}W_{m}^{\text{III},\alpha
}(x)dx=\frac{n\,\Gamma(n-m+\alpha+1)}{(n-m-1)!}\text{ for }n\geq m+1,\text{ as
required.}%
\]

To prove the norm formula in (\ref{eq:Type3norm}) for $n=0$ , we first
establish the following identity:
\begin{equation}
\int W_{m}^{\text{III},\alpha}(x)dx+m\int W_{m-1}^{\text{III},\alpha
-1}(x)dx=-\frac{x^{\alpha}e^{-x}}{L_{m}^{-\alpha-1}(-x)L_{m-1}^{-\alpha}(-x)}.
\label{WIII recursion}%
\end{equation}
Let $\psi(x)=L_{m}^{-\alpha-1}(-x)$ so that $\psi^{\prime}(x)=L_{m-1}%
^{-\alpha}(-x)$ and $\psi^{\prime\prime}(x)=L_{m-2}^{-\alpha+1}(-x).$ Now
$\psi(x)$ is a solution of the Laguerre differential equation%
\[
x\psi^{\prime\prime}(x)+(x-\alpha)\psi^{\prime}(x)-m\psi(x)=0.
\]
Divide this differential equation by $\psi(x)(\psi^{\prime}(x))^{2}$ and
rearrange to obtain%
\[
\frac{m}{(\psi^{\prime}(x))^{2}}+\frac{x}{\psi(x)}\left(  \frac{1}%
{\psi^{\prime}(x)}\right)  ^{\prime}+\frac{\alpha-x}{\psi(x)\psi^{\prime}%
(x)}=0.
\]
Multiplying by $x^{\alpha-1}e^{-x}$ yields%
\begin{equation}
\frac{mx^{\alpha-1}e^{-x}}{(\psi^{\prime}(x))^{2}}+\frac{x^{\alpha}e^{-x}%
}{\psi(x)}\left(  \frac{1}{\psi^{\prime}(x)}\right)  ^{\prime}+\frac{\left(
x^{\alpha}e^{-x}\right)  ^{\prime}}{\psi(x)\psi^{\prime}(x)}=0.
\label{Key Norm Identity}%
\end{equation}
Since%
\[
\frac{x^{\alpha}e^{-x}}{(\psi(x))^{2}}+\frac{x^{\alpha}e^{-x}}{\psi^{\prime
}(x)}\left(  \frac{1}{\psi(x)}\right)  ^{\prime}\equiv0,
\]
we see that \eqref{Key Norm Identity} can be rewritten as%
\begin{equation}
\frac{x^{\alpha}e^{-x}}{(\psi(x))^{2}}+\frac{mx^{\alpha-1}e^{-x}}%
{(\psi^{\prime}(x))^{2}}=-\left(  \frac{x^{\alpha}e^{-x}}{\psi^{\prime}%
(x)}\left(  \frac{1}{\psi(x)}\right)  ^{\prime}+\frac{x^{\alpha}e^{-x}}%
{\psi(x)}\left(  \frac{1}{\psi^{\prime}(x)}\right)  ^{\prime}+\frac{\left(
x^{\alpha}e^{-x}\right)  ^{\prime}}{\psi(x)\psi^{\prime}(x)}\right)  .
\label{Key Norm Identity 2}%
\end{equation}
From the product rule for derivatives, notice that
\begin{align}
\frac{x^{\alpha}e^{-x}}{\psi^{\prime}(x)}\left(  \frac{1}{\psi(x)}\right)
^{\prime}+\frac{x^{\alpha}e^{-x}}{\psi(x)}\left(  \frac{1}{\psi^{\prime}%
(x)}\right)  ^{\prime}+\frac{\left(  x^{\alpha}e^{-x}\right)  ^{\prime}}%
{\psi(x)\psi^{\prime}(x)}  &  =\left(  \left(  \frac{1}{\psi(x)}\right)
\cdot\left(  \frac{1}{\psi^{\prime}(x)}\right)  \cdot x^{\alpha}e^{-x}\right)
^{\prime}\label{Product Rule}\\
&  =\left(  \frac{x^{\alpha}e^{-x}}{L_{m}^{-\alpha-1}(-x)L_{m-1}^{-\alpha
}(-x)}\right)  ^{\prime}.\nonumber
\end{align}
Substituting \eqref{Product Rule} into \eqref{Key Norm Identity 2} and using
the definition of $W_{m}^{\text{III},\alpha}(x),$ we obtain%
\[
W_{m}^{\text{III},\alpha}(x)+mW_{m-1}^{\text{III},\alpha-1}(x)=-\left(
\frac{x^{\alpha}e^{-x}}{L_{m}^{-\alpha-1}(-x)L_{m-1}^{-\alpha}(-x)}\right)
^{\prime};
\]
integrating this expression now yields \eqref{WIII recursion}. Applying this
relation inductively yields
\begin{equation}
\int W_{m}^{\text{III},\alpha}(x)dx=m!(-1)^{m}\int e^{-x}x^{\alpha-m}%
dx-\sum_{j=0}^{m-1}(-1)^{j}\binom{m}{j}j!\frac{x^{\alpha-j}e^{-x}}%
{L_{m-j-1}^{-\alpha+j}(-x)L_{m-j}^{-\alpha+j-1}(-x)}. \label{eq:W3reduction}%
\end{equation}
For $r>0$, let $C_{r}=C_{1}+C_{2}+C_{3}$ denote the contour given by the ray
\[
C_{1}=\{x-ir:0\leq x<\infty\}
\]
oriented from right to left, by the left-side semi-circle
\[
C_{2}=\{re^{it}\colon\frac{\pi}{2}\leq t\leq\frac{3\pi}{2}\}
\]
oriented clockwise, and by the ray
\[
C_{3}=\{x+ir:0\leq x<\infty\}
\]
oriented from left to right. The zeros of $L_{m}^{-\alpha-1}(-x)$ are all
negative, and so by taking $r$ sufficiently small, the contour $C_{r}$ can be
made to not include these negative roots. Observe that the integrand
denominator is the square of a polynomial with simple roots. Hence the
residues of the integrand at the roots of $L_{m}^{-\alpha-1}(-x)$ vanish,
which means that it suffices to impose the condition that $r>0$ and that
$L_{m}^{-\alpha-1}(r)\neq0$. With this assumption,
\[
\int_{0}^{\infty}W_{m}^{\text{III},\alpha}(x)dx=\frac{1}{1-e^{2\pi i\alpha}%
}\int_{C_{r}}\frac{(-z)^{\alpha}e^{\pi i\alpha}\,e^{-z}}{\left(
L_{m}^{-\alpha-1}(-z)\right)  ^{2}}dz
\]
where $(-z)^{\alpha}$ denotes the principal branch of the power function. By
deforming $C_{r}$ we can rewrite the latter as a Mellin-Barnes integral,
namely
\begin{equation}
\int_{0}^{\infty}W_{m}^{\text{III},\alpha}(x)dx=\frac{i}{2\sin(\pi\alpha)}%
\int_{-r-i\infty}^{-r+i\infty}\frac{(-z)^{\alpha}e^{-z}}{\left(
L_{m}^{-\alpha-1}(-z)\right)  ^{2}}dz \label{eq:MBintW3}%
\end{equation}
Applying the same procedure to the usual integral representation of the
$\Gamma$-function gives
\begin{equation}
\frac{i}{2\sin(\pi a)}\int_{-r-i\infty}^{-r+i\infty}(-z)^{a}e^{-z}%
dz=\Gamma(1+a) \label{eq:MBintGamma}%
\end{equation}
valid for \emph{all} non-integral values of $a$ and all $r>0$. Applying
\eqref{eq:MBintW3} and \eqref{eq:MBintGamma} with $a=\alpha-m$ to
\eqref{eq:W3reduction} gives
\[
\int_{0}^{\infty}W_{m}^{\text{III},\alpha}(x)dx=m!(-1)^{m}\Gamma
(1+\alpha-m)=\frac{\Gamma(\alpha+1)\Gamma(-\alpha)m!}{\Gamma(m-\alpha)}.
\]
This completes the proof of the theorem.
\end{proof}
\end{theorem}

\subsection{The Completeness of the Type III Exceptional $X_{m}$-Laguerre
Polynomials\label{Completeness of Type III Laguerre Polynomials}}

In preparation for the proof of completeness of $\left\{  L_{m,n}%
^{\text{III},\alpha}\mid n=0,m+1,m+2,m+3,\ldots\right\}  $, we remind the
reader that the set $\mathcal{P}$ denotes the vector space of all polynomials
with complex coefficients in the real variable $x$ and, for $n\in
\mathbb{N}_{0},$ let $\mathcal{P}_{n}$ denote the vector space of all
$p\in\mathcal{P}$ with degree $\leq n.$ The following lemma is critical for
our argument; a proof can be found in \cite[Lemma 3, p. 416]{KMUG3}.

\begin{lemma}
\label{polynomial x P dense result}Suppose $\eta(x)$ is a polynomial such that
$\eta(x)\neq0$ for all $x\geq0.$ Then, for $\alpha>-1,$ the subspace%
\[
\eta\mathcal{P}:=\{\eta(x)p(x)\mid p\in\mathcal{P}\}
\]
is dense in $L^{2}((0,\infty);x^{\alpha}e^{-x}).$
\end{lemma}

We now prove the following completeness result.

\begin{theorem}
\label{Completeness} The set of Type $\text{III}$ exceptional $X_{m}$-Laguerre
polynomials%
\[
\left\{  L_{m,n}^{\textnormal{III},\alpha}\mid n=0,m+1,m+2,m+3,\ldots\right\}
\]
forms a complete orthogonal set of polynomials in the Hilbert space
$L^{2}((0,\infty);W_{m}^{\textnormal{III},\alpha}).$

\begin{proof}
The proof that we give of completeness is similar to the proof that we give of
completeness in Theorem \ref{Self-Adjoint Operator S Theorem}; we give the
full proof since some essential ingredients of the proof below are different.
\newline Let $\epsilon>0$ and let $f\in L^{2}((0,\infty);W_{m}^{\text{III}%
,\alpha}).$ Define%
\[
\widetilde{f}(x):=\frac{f(x)}{L_{m}^{-\alpha-1}(-x)};
\]
clearly%
\[
\left\Vert f\right\Vert _{m}^{\text{III},\alpha}=\left\Vert \widetilde{f}%
\right\Vert ^{\alpha},
\]
where $\left\Vert \cdot\right\Vert ^{\alpha}$ denotes the norm in
$L^{2}((0,\infty);x^{\alpha}e^{-x})$. Hence $\widetilde{f}\in L^{2}%
((0,\infty);x^{\alpha}e^{-x}).$ From Lemma \ref{polynomial x P dense result},
with%
\[
\eta(x)=L_{m}^{-\alpha-1}(-x),
\]
there exists $p\in\mathcal{P},$ say with \textrm{deg}$(p)=n,$ such that%
\[
\left\Vert \widetilde{f}(x)-L_{m}^{-\alpha-1}(-x)p(x)\right\Vert ^{\alpha
}<\epsilon.
\]
Hence it follows that%
\begin{align*}
\epsilon^{2}  &  >\int_{0}^{\infty}\left\vert \frac{f(x)}{L_{m}^{-\alpha
-1}(-x)}-L_{m}^{-\alpha-1}(-x)p(x)\right\vert ^{2}x^{\alpha}e^{-x}dx\\
&  =\int_{0}^{\infty}\left\vert \frac{f(x)-(L_{m}^{-\alpha-1}(-x))^{2}%
p(x)}{L_{m}^{-\alpha-1}(-x)}\right\vert ^{2}x^{\alpha}e^{-x}dx\\
&  =\int_{0}^{\infty}\left\vert f(x)-(L_{m}^{-\alpha-1}(-x))^{2}%
p(x)\right\vert ^{2}W_{m}^{\text{III},\alpha}(x)dx\\
&  =\left(  \left\Vert f(x)-(L_{m}^{-\alpha-1}(-x))^{2}p(x)\right\Vert
_{m}^{\text{III},\alpha}\right)  ^{2}.
\end{align*}
Notice that
\begin{equation}
(L_{m}^{-\alpha-1}(-x))^{2}p\in\mathcal{F}_{n+2m}, \label{F_(n+2m) condition}%
\end{equation}
where%
\[
\mathcal{F}_{n+2m}:=\left\{  P\in\mathcal{P}_{n+2m}\mid P^{\prime}%
(-x_{j})=0\text{ }(j=1,2,\ldots,m)\right\}  ,
\]
and where $\{x_{j}\}_{j=1}^{m}\subset(0,\infty)$ are the simple roots of the
Laguerre polynomial $L_{m}^{-\alpha-1}(x).$ We now show that $(L_{m}%
^{-\alpha-1}(-x))^{2}p(x)\in\mathcal{E}_{n+2m},$ where%
\[
\mathcal{E}_{n+2m}:=\mathrm{span}\left\{  L_{m,j}^{\text{III},\alpha}\mid
j=0,m+1,m+2,\ldots,n+2m\right\}  \text{;}%
\]
this will complete the proof of the theorem. Note that%
\begin{equation}
\text{\textrm{dim}}\left(  \mathcal{E}_{n+2m}\right)  =\text{\textrm{dim}%
}\left(  \mathcal{F}_{n+2m}\right)  =m+n+1. \label{dim E_(n+2m)}%
\end{equation}
Since the span of eigenfunctions of an operator is an invariant subspace of
that operator, we see that%
\[
\ell_{m}^{\text{III},\alpha}[\mathcal{E}_{n+2m}]\subset\mathcal{E}_{n+2m}.
\]
In particular, if $P\in\mathcal{E}_{n+2m},$ then
\begin{align*}
\ell_{m}^{\text{III},\alpha}[P](x)  &  =-xP^{\prime\prime}(x)+\left(
-1-\alpha+x+2x\frac{(L_{m}^{-\alpha-1}(-x))^{\prime}}{L_{m}^{-\alpha-1}%
(-x)}\right)  P^{\prime}(x)+(-m+\alpha)P(x)\\
&  \in\mathcal{P}.
\end{align*}
Consequently, the term%
\[
\frac{2x(L_{m}^{-\alpha-1}(-x))^{\prime}}{L_{m}^{-\alpha-1}(-x)}P^{\prime}(x)
\]
must be a polynomial. Since $L_{m}^{-\alpha-1}(x)$ is a classical Laguerre
polynomial, its roots $\{x_{j}\}_{j=1}^{m}$ $\subset(0,\infty)$ are simple so
it follows that $P^{\prime}(-x_{j})=0$ for $j=1,2,\ldots,m.$ Thus,%
\[
\mathcal{E}_{n+2m}\subset\mathcal{F}_{n+2m}.
\]
From $($\ref{dim E_(n+2m)}$)$ we see, in fact, that%
\[
\mathcal{E}_{n+2m}=\mathcal{F}_{n+2m}%
\]
From $($\ref{F_(n+2m) condition}$),$ it follows that
\[
(L_{m}^{-\alpha-1}(-x))^{2}p(x)\in\mathcal{E}_{n+2m},
\]
and this completes the proof.
\end{proof}
\end{theorem}

\subsection{Location of Roots and Root Asymptotics of the Type III Exceptional
$X_{m}$-Laguerre Polynomials\label{Properties of Roots}}

The following theorem gives us exact location of the $m+k$ (real) roots of the
Type III exceptional $X_{m}$-Laguerre polynomial $L_{m,m+k}^{\text{III}%
,\alpha}(x).$

\begin{theorem}
\label{Theorem Interlacing} Suppose $m,k\in\mathbb{N}$ and $-1<\alpha<0$. The
Type $\text{III}$ exceptional $X_{m}$-Laguerre polynomial $L_{m,m+k}%
^{\textnormal{III},\alpha}(x)$ has $k$ positive roots which interlace the roots of
$L_{k-1}^{\alpha+1}(x)$ and $m$ negative roots which interlace the roots of
$L_{m}^{-\alpha-1}(-x)$. More precisely, let $\{x_{k-1,i}^{\alpha+1}%
\}_{i=1}^{k-1}\subset(0,\infty)$ denote the $($simple$)$ roots of the Laguerre
polynomial $L_{k-1}^{\alpha+1}(x)$, and let $\{z_{m,i}^{-\alpha-1}\}_{i=1}%
^{m}\subset(-\infty,0)$ denote the $($simple$)$ roots of the Laguerre
polynomial $L_{m}^{-\alpha-1}(-x),$ with both sets ordered as follows:%
\[
z_{m,m}^{-\alpha-1}<z_{m,m-1}^{-\alpha-1}<\ldots<z_{m,1}^{-\alpha
-1}<0<x_{k-1,1}^{\alpha+1}<x_{k-1,2}^{\alpha+1}<\ldots<x_{k-1,k-1}^{\alpha
+1}.
\]
Then

\begin{enumerate}
\item[$(a)$] each of the $k$ intervals
\[
(0,x_{k-1,1}^{\alpha+1}),\text{ }(x_{k-1,1}^{\alpha+1},x_{k-1,2}^{\alpha
+1}),\ldots,\text{ }(x_{k-1,k-2}^{\alpha+1},x_{k-1,k-1}^{\alpha+1}),\text{
}(x_{k-1,k-1}^{\alpha+1},\infty)
\]
contains exactly one root of $L_{m,m+k}^{\textnormal{III},\alpha};$

\item[$(b)$] each of the $m$ intervals
\[
(-\infty,z_{m,m}^{-\alpha-1}),\text{ }(z_{m,m}^{-\alpha-1},z_{m,m-1}%
^{-\alpha-1}),\ldots,(z_{m,2}^{-\alpha-1},z_{m,1}^{-\alpha-1})
\]
contains exactly one root of $L_{m,m+k}^{\textnormal{III},\alpha}$.
\end{enumerate}

\begin{proof}
The key identity in establishing both $($a$)$ and $($b$)$ is the identity
given in $($\ref{Eq Alternative}$)$, namely%
\begin{align*}
L_{m,m+k}^{\text{III},\alpha}(x)  &  =(k+\alpha)L_{k-2}^{\alpha+1}%
(x)L_{m}^{-\alpha-1}(-x)+(m+1)L_{k-1}^{\alpha+1}(x)L_{m+1}^{-\alpha-1}(-x)\\
&  \quad-(m+k)L_{k-1}^{\alpha+1}(x)L_{m}^{-\alpha-1}(-x).
\end{align*}
We first prove part $($a$).$ Letting $x=x_{k-1,i}^{\alpha+1}$ $(i=1,2,\ldots
,k-1)$ yields
\begin{equation}
L_{m,m+k}^{\text{III},\alpha}(x_{k-1,i}^{\alpha+1})=(k+\alpha)L_{k-2}%
^{\alpha+1}(x_{k-1,i}^{\alpha+1})L_{m}^{-\alpha-1}(-x_{k-1,i}^{\alpha+1}).
\label{Location 1}%
\end{equation}
Since $L_{m}^{-\alpha-1}(-x)$ has no roots in $(0,\infty)$ and $L_{m}%
^{-\alpha-1}(0)>0,$ we see that
\[
L_{m}^{-\alpha-1}(-x_{k-1,i}^{\alpha+1})>0\quad(i=1,2,\ldots,k-1).
\]
Furthermore, from the classical theory, the roots of $L_{k-1}^{\alpha+1}(x)$
and $L_{k-2}^{\alpha+1}(x)$ interlace and since $L_{k-2}^{\alpha+1}(0)>0$, we
see that
\[
L_{k-2}^{\alpha+1}(x_{k-1,1}^{\alpha+1})>0.
\]
Hence, from $($\ref{Location 1}$),$ we deduce that
\begin{equation}
\text{sgn}(L_{m,m+k}^{\text{III},\alpha}(x_{k-1,i}^{\alpha+1}))=\text{sgn}%
(L_{k-2}^{\alpha+1}(x_{k-1,i}^{\alpha+1}))=(-1)^{i+1}\quad(i=1,\ldots,k-1).
\label{Type III Alternating Signs}%
\end{equation}
It follows that $L_{m,m+k}^{\text{III,}\alpha}(x)$ has a root in each of the
$k-2$ intervals%
\[
(x_{k-1,1}^{\alpha+1},x_{k-1,2}^{\alpha+1}),(x_{k-1,2}^{\alpha+1}%
,x_{k-1,3}^{\alpha+1}),\ldots,\text{ }(x_{k-1,k-2}^{\alpha+1},x_{k-1,k-1}%
^{\alpha+1}).
\]
From $($\ref{Negativity of Type III at 0}$),$ $L_{m,m+k}^{\text{III},\alpha
}(0)<0;$ hence, from $($\ref{Type III Alternating Signs}$)$ with $i=1,$ we see
that there is another root of $L_{m,m+k}^{\text{III},\alpha}(x)$ in the
interval $(0,x_{k-1,1}^{\alpha+1}).$ Lastly, from $($\ref{Lemma 2 Claim}$),$
we see that $x=x_{k-1,k-1}^{\alpha+1}$ is the right-most extreme point of
$L_{m,m+k}^{\text{III},\alpha}(x).$ Regardless of whether $x=x_{k-1,k-1}%
^{\alpha+1}$ is a relative maximum or relative minimum point, the graph of
$y=L_{m,m+k}^{\text{III},\alpha}(x)$ necessarily must cross the $x$-axis once
more at a point $x^{\ast}>x_{k-1,k-1}^{\alpha+1}.$ Hence we see that
$L_{m,m+k}^{\text{III},\alpha}(x)$ has a root in the interval $(x_{k-1,k-1}%
^{\alpha+1},\infty).$ Summarizing, we have shown that $L_{m,m+k}%
^{\text{III},\alpha}(x)$ has $k$ distinct, positive roots. \medskip\newline
The proof of $($b$)$ is similar. In this case, from $($\ref{Eq Alternative}%
$),$ we see that%
\begin{equation}
L_{m,m+k}^{\text{III},\alpha}(z_{m,i}^{-\alpha-1})=(m+1)L_{k-1}^{\alpha
+1}(z_{m,i}^{-\alpha-1})L_{m+1}^{-\alpha-1}(-z_{m,i}^{-\alpha-1}).
\label{Location 2}%
\end{equation}
Since $L_{k-1}^{\alpha+1}(x)$ has $k-1$ positive roots and $L_{k-1}^{\alpha
+1}(0)>0,$ we see that%
\begin{equation}
L_{k-1}^{\alpha+1}(z_{m,i}^{-\alpha-1})>0\quad(i=1,2,\ldots,m).
\label{Positivity at z's}%
\end{equation}
Moreover, since $L_{m+1}^{-\alpha-1}(0)>0,$ it follows from the interlacing
property of the roots of $L_{m}^{-\alpha-1}(x)$ and $L_{m+1}^{-\alpha-1}(x)$
that%
\[
L_{m+1}^{-\alpha-1}(-z_{m,i}^{-\alpha-1})=(-1)^{i}\quad(i=1,2,\ldots,m).
\]
Hence, from $($\ref{Location 2}$)$ and $($\ref{Positivity at z's}$),$ we see
that
\[
\text{sgn}(L_{m,m+k}^{\text{III},\alpha}(z_{m,i}^{-\alpha-1}))=\text{sgn}%
(L_{m+1}^{-\alpha-1}(-z_{m,i}^{-\alpha-1}))=(-1)^{i}\quad(i=1,\ldots,m).
\]
This implies that each of the $m-1$ intervals%
\[
(z_{m,m}^{-\alpha-1},z_{m,m-1}^{-\alpha-1}),\ldots,(z_{m,2}^{-\alpha
-1},z_{m,1}^{-\alpha-1})
\]
contains a root of $L_{m,m+k}^{\text{III},\alpha}(x).$ We claim that there is
an additional root of $L_{m,m+k}^{\text{III},\alpha}(x)$ in the interval
$(-\infty,z_{m,m}^{-\alpha-1}).$ Indeed, from $($\ref{Lemma 2 Claim}$),$
$x=z_{m,m}^{-\alpha-1}$ is the left-most extreme point of $L_{m,m+k}%
^{\text{III},\alpha}(x)$ so, as in part $($a$)$, there must be another root of
$L_{m,m+k}^{\text{III},\alpha}(x)$ at a point $z^{\ast}<$ $z_{m,m}^{-\alpha
-1}.$ This completes the proof that $L_{m,m+k}^{\text{III},\alpha}(x)$ has $m$
roots in $(-\infty,0).$ Combining this fact with part $($a$),$ we have found
all $m+k$ roots of $L_{m,m+k}^{\text{III},\alpha}(x)$ and this completes the
proof of the theorem.
\end{proof}
\end{theorem}

\begin{remark}
When $k=1,$ $L_{m,m+1}^{\textnormal{III},\alpha}(x)$ has one positive root; the
exact location of this root cannot be specifically identified. When $m=1,$ the
Laguerre polynomial $L_{1}^{-\alpha-1}(-x)$ has a unique root $z_{1,1}%
^{-\alpha-1}<0.$ In this case, the above theorem indicates there is a unique
root of $L_{1,k+1}^{\textnormal{III},\alpha}(x)$ in the interval $(-\infty
,z_{1,1}^{-\alpha-1}).$
\end{remark}

We call the $m$ negative roots of $L_{m,m+k}^{\text{III,}\alpha}(x)$ above the
`exceptional' roots of $L_{m,m+k}^{\text{III},\alpha}(x).$ We now discuss the
asymptotic behavior of the roots as $k\rightarrow\infty$.

\begin{theorem}
\label{Theorem Asymptotic} As $k\rightarrow\infty$:

\begin{itemize}
\item[$(a)$] The exceptional roots of $L_{m,m+k}^{\textnormal{III},\alpha}(x)$
converge to the roots of $L_{m}^{-\alpha-1}(-x).$

\item[$(b)$] The first positive root of $L_{m,m+k}^{\textnormal{III},\alpha}(x)$
tends to zero.
\end{itemize}

\begin{proof}
Recall, from $($\ref{Type III Laguerre rep}$),$ that
\begin{equation}
L_{m,m+k}^{\text{III},\alpha}(x)=xL_{k-2}^{\alpha+2}(x)L_{m}^{-\alpha
-1}(-x)+(m+1)L_{k-1}^{\alpha+1}(x)L_{m+1}^{-\alpha-2}(-x).
\label{Theorem Asymptotic 1}%
\end{equation}
Now, the outer ratio asymptotics for the classical Laguerre polynomials give
\[
\frac{L_{k-1}^{\alpha+1}(x)}{L_{k-2}^{\alpha+2}(x)}\,\,\simeq\left(  -\frac
{x}{k}\right)  ^{1/2}+O\left(  \frac{1}{k}\right)  ,\quad k\rightarrow
\infty\quad
\]
with convergence uniform on compact sets that avoid the positive real axis.
Therefore, dividing the identity in $($\ref{Theorem Asymptotic 1}$)$ by
$L_{k-2}^{\alpha+2}(x)$ and taking the limit as $k\rightarrow\infty,$ we
obtain
\[
\frac{L_{m,m+k}^{\text{III},\alpha}(x)}{L_{k-2}^{\alpha+2}(x)}%
\,\,\overset{k\rightarrow\infty}{\longrightarrow}\,\,xL_{m}^{-\alpha-1}(-x).
\]
Part $(a)$ follows by Hurwitz's theorem $($\cite[Theorem 1.91.3]{Szego}$)$ and
Theorem \ref{Theorem Interlacing}.\newline The extra $x$ factor implies Part
$(b)$ as follows. At $x=0$ we have%
\[
\frac{L_{k-1}^{\alpha+1}(0)}{L_{k-2}^{\alpha+1}(0)}=\frac{k+\alpha}%
{k-1}\overset{k\rightarrow\infty}{\longrightarrow}1.
\]
So the ratio asymptotics hold at $x=0.$ And from $($%
\ref{Simplest Representation}$)$, it follows that $L_{m,m+k}^{\text{III}%
,\alpha}(0)\rightarrow0$ as $k\rightarrow\infty.$
\end{proof}
\end{theorem}

\subsection{Type $\text{III}$ Exceptional $X_{m}$-Laguerre Spectral
Analysis\label{Type III Spectral Analysis}}

The maximal domain associated with the differential expression $\ell
_{m}^{\text{III,}\alpha}[\cdot],$ given in either (\ref{Type III Expression}%
)\ or (\ref{Type III Symmetric Form}), in the Hilbert space $L^{2}\left(
(0,\infty);W_{m}^{\text{III},\alpha}\right)  $ is defined to be
\[
\Delta_{m}^{\text{III,}\alpha}:=\left\{  f:(0,\infty)\rightarrow\mathbb{C}\mid
f,f^{\prime}\in AC_{\text{\textrm{loc}}}(0,\infty);f,\ell_{m}^{\text{III,}%
\alpha}[f]\in L^{2}((0,\infty);W_{m}^{\text{III},\alpha})\right\}  .
\]
The associated maximal operator
\[
T_{1,m}^{\text{III,}\alpha}:\mathcal{D}(T_{1,m}^{\text{III,}\alpha})\subset
L^{2}((0,\infty);W_{m}^{\text{III},\alpha})\rightarrow L^{2}((0,\infty
);W_{m}^{\text{III},\alpha})\,,
\]
is defined to be
\begin{align*}
T_{1,m}^{\text{III,}\alpha}f  &  =\ell_{m}^{\text{III,}\alpha}[f]\\
f\in\mathcal{D}(T_{1,m}^{\text{III,}\alpha}):  &  =\Delta_{m}^{\text{III,}%
\alpha}.
\end{align*}
For $f,g\in\Delta_{m}^{\text{III,}\alpha}$, Green's formula may be written as
\[
\int_{0}^{\infty}\ell_{m}^{\text{III,}\alpha}[f](x)\overline{g}(x)W_{m}%
^{\text{III},\alpha}(x)dx=[f,g]_{m}^{\text{III,}\alpha}(x)\mid_{x=0}%
^{x=\infty}+\int_{0}^{\infty}f(x)\ell_{m}^{\text{III,}\alpha}[\overline
{g}](x)W_{m}^{\text{III},\alpha}(x)dx\,,
\]
where $[\cdot,\cdot]_{m}^{\text{III},\alpha}(\cdot)$ is the sesquilinear form
defined by
\[
\left[  f,g\right]  _{m}^{\text{III,}\alpha}(x):=\frac{x^{\alpha+1}e^{-x}%
}{(L_{m}^{-\alpha-1}(-x))^{2}}(f(x)\overline{g}^{\prime}(x)-f^{\prime
}(x)\overline{g}(x))\quad(0<x<\infty)
\]
and where
\begin{align*}
\left[  f,g\right]  _{m}^{\text{III},\alpha}(x)\mid_{x=0}^{x=\infty}:  &
=\left[  f,g\right]  _{m}^{\text{III},\alpha}(\infty)-\left[  f,g\right]
_{m}^{\text{III,}\alpha}(0)\,\\
&  :=\lim_{x\rightarrow\infty}\left[  f,g\right]  _{m}^{\text{III},\alpha
}(x)-\lim_{x\rightarrow0^{+}}\left[  f,g\right]  _{m}^{\text{III},\alpha}(x)
\end{align*}

The adjoint of the maximal operator in $L^{2}\left(  (0,\infty);W_{m}%
^{\text{III},\alpha}\right)  $ is the minimal operator $T_{0,m}^{\text{III}%
,\alpha},$ defined by%
\begin{align*}
T_{0,m}^{\text{III},\alpha}f  &  =\ell_{m}^{\text{III,}\alpha}[f]\\
f\in\mathcal{D}(T_{0,m}^{\text{III},\alpha}):  &  =\{f\in\Delta_{m}%
^{\text{III},\alpha}\mid\lbrack f,g]_{m}^{\text{III},\alpha}\mid
_{x=0}^{x=\infty}=0\text{ for all }g\in\Delta_{m}^{\text{III},\alpha}\}.
\end{align*}

Both endpoints $x=0$ and $x=\infty$ are singular points of $\ell
_{m}^{\text{III,}\alpha}[\cdot].$ In fact, $x=0$ is a regular singular
endpoint in the sense of Frobenius and $x=\infty$ is an irregular singular
endpoint. The associated Frobenius indicial equation at $x=0$ is
$r(r+\alpha)=0.$ Consequently, two linearly independent solutions of $\ell
_{m}^{\text{III,}\alpha}[y]=0$ will behave asymptotically like
\[
z_{1}(x):=1\quad\text{and}\quad z_{2}(x):=x^{-\alpha}%
\]
near $x=0.$ Since $-1<\alpha<0,$ it is clear that both solutions are in
$L^{2}((0,\infty);W_{m}^{\text{III},\alpha});$ in other words, $\ell
_{m}^{\text{III,}\alpha}[\cdot]$ is in the limit-circle case at $x=0.$

For the analysis at the irregular singular endpoint, $x=\infty$, we obtain two
linearly independent solutions using the standard reduction of order method.
Solving the differential equation $\ell_{m}^{\text{III,}\alpha}[y](x)=0$ we
have a basis of solutions $\left\{  y_{1}(x),y_{2}(x)\right\}  $, where
\[
y_{1}(x)=1\in L^{2}((0,\infty);W_{m}^{\text{III},\alpha})
\]
and
\[
y_{2}(x)=\int_{a}^{x}\frac{e^{t}(L_{m}^{-\alpha-1}(t))^{2}}{t^{\alpha+1}%
}\,dt\,\quad(a>0\text{ is arbitrary).}%
\]
Mimicking the proof of the Type $\text{I}$ case in Section
\ref{Type I Spectral Analysis}, we find that $y_{2}\notin L^{2}((1,\infty
);W_{m}^{\text{III},\alpha})$. Consequently, we obtain the following result on
the deficiency indices of $T_{0,m}^{\text{III},\alpha}.$

\begin{theorem}
Let $T_{0,m}^{\textnormal{III},\alpha}$ be the minimal operator in $L^{2}%
((0,\infty);W_{m}^{\textnormal{III},\alpha})$ generated by the Type$\ $III
exceptional $X_{m}$-Laguerre differential expression $\ell_{m}%
^{\textnormal{III,}\alpha}[\cdot]$. For $-1<\alpha<0$, the deficiency index of
$T_{0,m}^{\textnormal{III},\alpha}$ is $(1,1)$.
\end{theorem}

For $-1<\alpha<0$, we must impose one boundary condition at $x=0$ in order to
obtain a self-adjoint extension of the minimal operator $T_{0,m}%
^{\text{III},\alpha}$. We seek to find that self-adjoint operator
$T_{m}^{\text{III},\alpha}$ which has the Type III polynomials $\left\{
L_{m,n}^{\text{III},\alpha}\mid n=0,m+1,m+2,m+3,\ldots\right\}  $ as eigenfunctions.

Note that $x^{-\alpha}$ $\in L^{2}((0,\infty);W_{m}^{\text{III},\alpha})$
since $-1<\alpha<0$. A calculation shows that
\[
\ell_{m}^{\text{III,}\alpha}[x^{-\alpha}]=\frac{-2\alpha x^{-\alpha}%
L_{m-1}^{-\alpha}(-x)}{L_{m}^{-\alpha-1}(-x)}-mx^{-\alpha}.
\]
Since%
\begin{align*}
&  \int_{0}^{\infty}\left\vert \frac{x^{-\alpha}L_{m-1}^{-\alpha}(-x)}%
{L_{m}^{-\alpha-1}(-x)}\right\vert ^{2}\frac{x^{\alpha}e^{-x}}{(L_{m}%
^{-\alpha-1}(-x))^{2}}dx\medskip\\
&  \leq\frac{1}{\left(  L_{m}^{-\alpha-1}(0)\right)  ^{4}}\int_{0}^{\infty
}\left(  L_{m-1}^{-\alpha}(-x)\right)  ^{2}x^{-\alpha}e^{-x}dx\medskip\\
&  <\infty,
\end{align*}
we see that $\ell_{m}^{\text{III,}\alpha}[x^{-\alpha}]\in L^{2}((0,\infty
);W_{m}^{\text{III},\alpha}).$ Consequently, $x^{-\alpha}\in\Delta
_{m}^{\text{III,}\alpha}$ for $-1<\alpha<0$. Moreover, the calculation
\[
\lbrack x^{-\alpha},1]_{m}^{\text{III,}\alpha}(0)=\alpha\lim_{x\rightarrow
0}\frac{e^{-x}}{(L_{m}^{-\alpha-1}(-x))^{2}}\neq0\,,
\]
proves that $1\notin\mathcal{D}(T_{0,m}^{\text{III,}\alpha}),$ the minimal
domain, and thus we can use the function $1$ as an appropriate Glazman
boundary function. For $f\in\Delta_{m}^{\text{III,}\alpha},$ further
calculations show that%
\[
0=[f,1]_{m}^{\text{III,}\alpha}(0)=\lim_{x\rightarrow0^{+}}x^{\alpha
+1}f^{\prime}(x)
\]
and%
\[
\lim_{x\rightarrow0^{+}}x^{\alpha+1}(L_{m,n}^{\text{III,}\alpha}(x))^{\prime
}=0.
\]
Summarizing, and using Theorem \ref{Completeness}, we obtain the following theorem.

\begin{theorem}
\label{Type III Self-Adjoint Operator}Suppose $-1<\alpha<0$. The operator
\[
T_{m}^{\textnormal{III},\alpha}:\mathcal{D}(T_{m}^{\textnormal{III},\alpha})\subset
L^{2}((0,\infty);W_{m}^{\textnormal{III},\alpha})\rightarrow L^{2}((0,\infty
);W_{m}^{\textnormal{III},\alpha})\,,
\]
defined by
\begin{align*}
T_{m}^{\textnormal{III},\alpha}f  &  =\ell_{m}^{\textnormal{III,}\alpha}[f]\\
f\in\mathcal{D}(T_{m}^{\textnormal{III},\alpha}):  &  =\{f\in\Delta_{m}%
^{\textnormal{III,}\alpha}\mid\lim_{x\rightarrow0^{+}}x^{\alpha+1}f^{\prime
}(x)=0\},
\end{align*}
is a self-adjoint extension of the minimal operator $T_{0,m}%
^{\textnormal{III,}\alpha}$ in $L^{2}((0,\infty);W_{m}^{\textnormal{III},\alpha})$
having the Type III exceptional $X_{m}$-Laguerre polynomials%
\[
\{L_{m,n}^{\textnormal{III},\alpha}\mid n=0,m+1,m+2,m+3,\ldots\}
\]
as a complete set of eigenfunctions. Moreover, the spectrum of $T_{m}%
^{\textnormal{III},\alpha}$ consists only of eigenvalues and is given explicitly
by
\[
\sigma(T_{m}^{\textnormal{III},\alpha})=\left\{  n-m+\alpha\mid
n=0,m+1,m+2,m+3,\ldots\right\}  \,\,.
\]

\end{theorem}

\section{Appendix}

The following is a list of a few Type $\text{III}$ exceptional $X_{m}%
$-Laguerre polynomials $\left\{  L_{m,n}^{\text{III,}\alpha}\right\}  $ for
various values of $m$ and $n$. A similar list of Type $\text{I}$ and
$\text{II}$ exceptional $X_{m}$-Laguerre polynomials can be found in
\cite{HoSasakiZeros2012}.

For $m=1$ we have:%

\begin{align*}
L_{1,0}^{\text{III},\alpha}(x)  &  =1\\
L_{1,2}^{\text{III},\alpha}(x)  &  =x^{2}-2\alpha x+\alpha(\alpha+1)\\
L_{1,3}^{\text{III},\alpha}(x)  &  =-x^{3}+3(\alpha+1)x^{2}-3\alpha
(\alpha+2)x+\alpha(\alpha+1)(\alpha+2)\\
L_{1,4}^{\text{III},\alpha}(x)  &  =\frac{1}{2}x^{4}-2(\alpha+2)x^{3}%
+(\alpha+3)(3\alpha+2)x^{2}-2\alpha(\alpha+2)(\alpha+3)x\\
&  \quad\quad\quad\quad+\frac{1}{2}\alpha(\alpha+1)(\alpha+2)(\alpha+3)\\
L_{1,5}^{\text{III},\alpha}(x)  &  =-\frac{1}{6}x^{5}+\frac{5}{6}%
(\alpha+3)x^{4}-\frac{5}{6}(\alpha+4)(2\alpha+3)x^{3}\\
&  \quad\quad\quad\quad+\frac{5}{6}(\alpha+3)(\alpha+4)(2\alpha+1)x^{2}%
-\frac{5}{6}\alpha(\alpha+2)(\alpha+3)(\alpha+4)x\\
&  \quad\quad\quad\quad+\frac{1}{6}\alpha(\alpha+1)(\alpha+2)(\alpha
+3)(\alpha+4).
\end{align*}

For $m=2$ we obtain:%

\begin{align*}
L_{2,0}^{\text{III},\alpha}(x) &  =1\\
L_{2,3}^{\text{III},\alpha}(x) &  =\frac{1}{2}x^{3}-\frac{3(\alpha-1)}{2}%
x^{2}+\frac{3\alpha(\alpha-1)}{2}x-\frac{\alpha(\alpha-1)(\alpha+1)}{2}\\
L_{2,4}^{\text{III},\alpha}(x) &  =-\frac{1}{2}x^{4}+2\alpha x^{3}%
-(\alpha-1)(3\alpha+4)x^{2}+2\alpha(\alpha-1)(\alpha+2)x\\
&  \quad\quad\quad\quad-\frac{\alpha(\alpha-1)(\alpha+1)(\alpha+2)}{2}\\
L_{2,5}^{\text{III},\alpha}(x) &  =\frac{1}{4}x^{5}-\frac{5(\alpha+1)}{4}%
x^{4}+\frac{5}{2}(\alpha^{2}+2\alpha-1)x^{3}\\
&  -\frac{5(\alpha-1)(\alpha+1)(\alpha+3)}{2}x^{2}+\frac{5\alpha
(\alpha-1)(\alpha+2)(\alpha+3)}{4}x\\
&  -\frac{\alpha(\alpha-1)(\alpha+1)(\alpha+2)(\alpha+3)}{4}%
\end{align*}

\begin{align*}
L_{2,6}^{\text{III},\alpha}(x)  &  =-\frac{1}{12}x^{6}+\frac{(\alpha+2)}%
{2}x^{5}-\frac{(5\alpha^{2}+19\alpha+6)}{4}x^{4}\\
&  \quad\quad\quad\quad+\frac{(\alpha+2)(\alpha+4)(5\alpha-3)}{3}x^{3}%
-\frac{(\alpha-1)(\alpha+3)(\alpha+4)(5\alpha+4)}{4}x^{2}\\
&  \quad\quad\quad\quad+\frac{\alpha(\alpha-1)(\alpha+2)(\alpha+3)(\alpha
+4)}{2}x\\
&  \quad\quad\quad\quad-\frac{\alpha(\alpha-1)(\alpha+1)(\alpha+2)(\alpha
+3)(\alpha+4)}{12}.
\end{align*}

For $m=3$:
\begin{align*}
L_{3,0}^{\text{III},\alpha}(x)  &  =1\\
L_{3,4}^{\text{III},\alpha}(x)  &  =\frac{1}{6}x^{4}-\frac{2(\alpha-2)}%
{3}x^{3}+(\alpha-1)(\alpha-2)x^{2}\\
&  \quad\quad\quad\quad-\frac{2\alpha(\alpha-1)(\alpha-2)}{3}x+\frac
{\alpha(\alpha-2)(\alpha-1)(\alpha+1)}{6}\\
L_{3,5}^{\text{III},\alpha}(x)  &  =-\frac{1}{6}x^{5}+\frac{5(\alpha-1)}%
{6}x^{4}-\frac{5(\alpha-2)(2\alpha+1)}{6}x^{3}\\
&  \quad\quad\quad\quad+\frac{5(\alpha-2)(\alpha-1)(2\alpha+3)}{6}x^{2}%
-\frac{5\alpha(\alpha-2)(\alpha-1)(\alpha+2)}{6}x\\
&  \quad\quad\quad\quad+\frac{\alpha(\alpha-2)(\alpha-1)(\alpha+1)(\alpha
+2)}{6}\\
L_{3,6}^{\text{III},\alpha}(x)  &  =\frac{1}{12}x^{6}-\frac{\alpha}{2}%
x^{5}+\frac{(5\alpha^{2}+\alpha-12)}{4}x^{4}-\frac{\alpha(\alpha
-2)(5\alpha+13)}{3}x^{3}\\
&  \quad\quad\quad\quad+\frac{(\alpha-2)(\alpha-1)(\alpha+3)(5\alpha+6)}%
{4}x^{2}\\
&  \quad\quad\quad\quad-\frac{\alpha(\alpha-2)(\alpha-1)(\alpha+2)(\alpha
+3)}{2}x\\
&  \quad\quad\quad\quad+\frac{\alpha(\alpha-2)(\alpha-1)(\alpha+1)(\alpha
+2)(\alpha+3)}{12}\\
L_{3,7}^{\text{III},\alpha}(x)  &  =-\frac{1}{36}x^{7}+\frac{7(\alpha+1)}%
{36}x^{6}-\frac{7(\alpha^{2}+2\alpha-2)}{12}x^{5}\\
&  \quad\quad\quad\quad+\frac{35(\alpha+1)(\alpha^{2}+2\alpha-6)}{36}%
x^{4}-\frac{7(\alpha-2)(\alpha+4)(5\alpha^{2}+10\alpha-3)}{36}x^{3}\\
&  \quad\quad\quad\quad+\frac{7(\alpha-2)(\alpha-1)(\alpha+1)(\alpha
+3)(\alpha+4)}{12}x^{2}\\
&  \quad\quad\quad\quad-\frac{7\alpha(\alpha-2)(\alpha-1)(\alpha
+2)(\alpha+3)(\alpha+4)}{36}x\\
&  \quad\quad\quad\quad+\frac{\alpha(\alpha-2)(\alpha-1)(\alpha+1)(\alpha
+2)(\alpha+3)(\alpha+4)}{36}.
\end{align*}

\end{document}